\documentclass[journal]{IEEEtran}

\usepackage{amssymb}
\usepackage{amsfonts}
\usepackage{latexsym, amssymb,amsmath, amsbsy,amsopn, amstext}
\usepackage{amsmath, amsbsy}
\usepackage{hyperref}
\usepackage[none]{hyphenat}
\usepackage{float}
\usepackage{ifpdf}
\usepackage{xcolor}
\usepackage{tikz}
\usetikzlibrary{arrows,shapes,snakes,shadows,positioning,automata,patterns}
\usetikzlibrary{trees,decorations.pathmorphing,decorations.markings}

\def\diag{diag}

\def\0{{\bf 0}}

\newtheorem{thm}{Theorem}
\newtheorem{lem}{Lemma}
\newtheorem{cor}{Corollary}

\newtheorem{rem}{Remark}
\newtheorem{alg}{Algorithm}
\newtheorem{assum}{Assumption}

\title{\LARGE \bf
Asynchronous distributed algorithms for seeking generalized Nash equilibria under full  and partial decision information
}

\author{Peng Yi  and Lacra Pavel
\thanks{This work was supported by NSERC Discovery Grant (261764). A preliminary  version of this work will be presented at 2018 European Control Conference.}
\thanks{P. Yi is with Department of Electrical and Systems Engineering, Washington University in St. Louis, USA,
 and L. Pavel is with Department of Electrical and Computer Engineering, University of Toronto, Canada.
        {\tt\small peng.yi@utoronto.edu,pavel@control.toronto.edu}}%
}

\begin{document}

\maketitle
\thispagestyle{empty}
\pagestyle{empty}

\begin{abstract}
We investigate asynchronous distributed algorithms with delayed information for seeking generalized Nash equilibrium over multi-agent networks. The considered game model has all players' local decisions coupled with a shared affine constraint. We assume each player can only access its local objective function, local constraint, and a local block matrix of the affine constraint.  We first give the algorithm  for the case when each agent is able to fully access all other players' decisions.
With the help of auxiliary edge variables and {\it edge Laplacian matrix}, each player can carry on its local iteration in an asynchronous manner, using only local data and possibly delayed neighbour information.
And then we investigate the case when the agents cannot know all other players' decisions, which is called a {\it partial decision case}.
We introduce a local estimation of the overall decisions for each agent in the partial decision case, and develop another asynchronous algorithm by incorporating consensus dynamics on the local estimations of the overall decisions.
 Since both algorithms do not need any centralized clock coordination, the algorithms fully exploit the local computation resource of each player, and remove the idle time due to waiting for the ``slowest" agent.
 Both  algorithms are developed by a preconditioned forward-backward operator splitting method, while  convergence is shown with the help of asynchronous fixed-point iterations under proper assumptions and fixed step-size choices.
 Numerical studies verify both algorithms' convergence and efficiency.
\end{abstract}

\begin{IEEEkeywords}
monotone game, generalized Nash equilibrium, distributed algorithm, asynchronous computation, multi-agents systems, operator splitting methods
\end{IEEEkeywords}

\IEEEpeerreviewmaketitle

\section{Introductions }\label{sec_introduction}

In many large-scale noncooperative multi-agent networks, the agents (players) need to make local decisions to minimize/maximize their individual cost/utility, for e.g., power allocation in cognitive radio networks, \cite{pang,scutari}, demand response and electric-vehicle charging management in smart grids, \cite{hu,zhuminghui,grammatico_2,lygeros2},
rate/power control over optical networks, \cite{pavel1,pavel2}, and opinion evolution over social networks, \cite{opion,lygeros1}.
Each player's local objective function to be optimized depends on other players' decisions, and
their feasible decision sets may also couple through shared/non-shared constraints.
Generalized Nash Equilibrium (GNE) is  the reasonable solution to noncooperative games with coupling constraints,
where  no player can  decrease/increase its  cost/utility by unilaterally changing its local decision to another feasible one.
Interested readers can refer to \cite{kanzow,faccinei2} for a review on GNE.

Recently, distributed NE/GNE computation methods over large-scale networks have received increasing  attention, see \cite{hu}-\cite{lygeros2} and \cite{peng4}-\cite{yipeng2}.
With distributed GNE seeking methods, each player keeps its local objective function, local feasible set and local constraints, and the data does not need to be transferred to a center.
Various game models and information structures have been considered in those works. \cite{zhuminghui} proposed a primal-dual approach for GNE seeking with diminishing step-sizes for pseudo-monotone games. \cite{grammatico_2,lygeros2} and \cite{grammatico_1} proposed  fixed-point iteration approaches for NE/GNE computation of aggregative games with/without affine coupling constraints, but also adopted a centralized coordinator.    \cite{shanbhag4} considered GNE seeking for monotone games based on a Tikhonov regularized primal-dual algorithm with a centralized multiplier coordinator.
\cite{yipeng} proposed a center-free  GNE seeking algorithm for strongly monotone games with fixed step-sizes.  \cite{yipeng2} proposed an ADMM approach for seeking GNE of  monotone games based on preconditioned proximal algorithms. The previously mentioned works assumed that each player is able to
know all other players' decisions or other players' decisions that its local objective function {\it directly} depends on.
We call these as  algorithms with  {\it full decision information}. On the other hand, \cite{shanbhag1}, \cite{pavel3},\cite{pavel4},\cite{pavel5} and \cite{pavel6} considered  gradient-based methods for NE/GNE seeking, and assumed that each player can only estimate all other players' decisions through  local communications.
\cite{hu,hu1,hu2} also considered NE seeking with a partial decision information, and combined  continuous-time consensus dynamics and gradient flow  to ensure that the local estimations will reach consensus on the NE.
\cite{shanbhag1} considered NE computation for aggregative games when the players can only get a local estimation for the aggregative variable. And
\cite{liangshu} combined  projected  continuous-time dynamics and finite-time consensus dynamics for GNE seeking of aggregative games  by introducing local estimators for aggregative variables.
\cite{wanglong} proposed a continuous-time dynamics for seeking GNE by incorporating a leader-follower consensus algorithm when each player only has an estimation of other players' decisions.
Since these algorithms do not require that each  agent  knows all other players' decisions, we call them  algorithms with  {\it partial decision information}.

However,  all the above works considered only {\it synchronous} GNE computation algorithms,
which need a global clock or a coordinator to ensure that all players have finished their current iterations before executing the next one.
This is quite computationally inefficient since all players have to wait for the ``slowest" player to finish before carrying on their local iterations.
Meanwhile, the centralized coordinator needs to guarantee that for their computations, all players utilize the most recent information at the same iteration index.
Hence, the synchronicity requirement limits their application and efficiency in large-scale networks. In fact, asynchronous computation has always been an appealing feature for algorithms over multi-agent networks. For example,  \cite{yinwotao} proposed a distributed optimization algorithm with asynchrony and delays, which was shown to enjoy a high computational efficiency compared to the synchronous ones. \cite{scutari} investigated asynchronous NE computation for monotone games, and \cite{jinlong} investigated asynchronous NE computation for stochastic games, both with proximal best-response algorithms.
Moreover, naively running the synchronous algorithms in an asynchronous manner will not always enjoy convergence properties.
However, the continuous-time algorithms in \cite{hu1,hu2,wanglong,liangshu} require synchronization of time and continuous communication among the players,
therefore, are not directly implementable with digital communications and  computations.
Hence, it is important to design and analyze discrete-time asynchronous GNE seeking algorithms with convergence guarantees.

Motivated by the above, the main contribution of this paper is to propose  novel \emph{asynchronous distributed GNE} seeking algorithms for noncooperative games with linear equality coupling constraints, under {\it full decision information} and {\it partial decision information}, respectively.
The equality constraint is motivated by various
task allocation problems when a global requirement should be exactly met by all players' decisions, like electric-vehicle charging management \cite{hu,grammatico_2}.
The algorithms are  motivated by a preconditioned forward-backward operator splitting method for finding zeros of a sum of monotone operators, which was also adopted in \cite{yipeng}. However, due to the sequential updating order of the algorithm in \cite{yipeng}, it is nontrivial to develop an asynchronous algorithm based on \cite{yipeng}.
Moreover, the algorithm in \cite{yipeng} requires full decision information.
In this paper, we first achieve asynchronous GNE  computation with {\it full decision information} by introducing   auxiliary {\it edge variables} and local multipliers and using {\it edge Laplacian} matrix,  unlike  \cite{yipeng}, where nodal variables and Laplacian are used. The proposed  asynchronous distributed algorithm can rely  on delayed information, thus further eliminates the need for coordinating information consistence in terms of iteration index.
Moreover, the proposed algorithm eliminates idle time, and each player can carry on its local iteration as soon as the delayed neighbour information is available. Assumptions similar to those in \cite{yipeng}, \cite{sayed}, \cite{grammatico_2} and \cite{lygeros2} are adopted for its convergence analysis.
Then we propose an asynchronous algorithm with {\it partial decision information},  where each agent  maintains a local estimation of the overall decisions.
To ensure that the agents' local estimations reach the same GNE, an additional consensus dynamics is incorporated into the algorithm. Assumptions similar to \cite{pavel5}, \cite{pavel6} and \cite{tatarenko} are adopted for the algorithm with {\it partial decision information}.
Both algorithms are shown to be developed by the preconditioned forward-backward operator splitting method with properly chosen operators.
The convergence analysis is given  based on \cite{yinwotao_arock} for asynchronous fixed-point iterations,  under bounded delay assumption and  a sufficient  {\it fixed step-size} choice. Compared with the preliminary conference paper \cite{yipeng3}, we provide the additional algorithm with {\it partial decision information} and the detailed analysis and convergence proofs.

The paper is organized as follows.
Section \ref{sec_notations_and preliminaries} gives the notations and preliminary background.
Section \ref{sec_game_formulation} formulates the noncooperative game.
Section \ref{sec_algorithm_and_limiting_poinit} presents our first asynchronous distributed GNE computation algorithm with {\it full decision information},
and Section \ref{sec_convergence_analysis} presents its development and convergence analysis.
Section \ref{sec_pdi} presents the algorithm and analysis for the case with {\it partial decision information}.
Section \ref{sec_numerical_studies} presents  numerical  studies for both algorithms, and
Section \ref{sec_concluding} draws the concluding remarks.

\section{notations and preliminaries}\label{sec_notations_and preliminaries}

In this section, we review the notations and preliminary notions in  monotone and averaged operators from \cite{combettes1}.

{\it Notations}:  In the following,
$\mathbf{R}^m$ ($\mathbf{R}^m_{+}$) denotes the $m-$dimesional (nonnegative) Euclidean space.
For a column vector $x \in \mathbf{R}^m$ (matrix $A\in \mathbf{R}^{m\times n}$),
$x^T$ ($A^T$)  denotes its transpose.
$x^Ty=\langle x, y\rangle$ denotes the inner product of $x,y$, and $||x||= \sqrt{x^Tx}$ denotes the norm induced by inner product $\langle\cdot,\cdot\rangle$.
$||x||^2_G$ denotes $\langle x, Gx\rangle$ for a symmetric positive definite matrix $G$,
and $||x||_G=\sqrt{\langle x, Gx\rangle}$ denotes $G-$induced norm.
Denote $\mathbf{1}_m=(1,...,1)^T \in \mathbf{R}^m$ and
$\mathbf{0}_m=(0,...,0)^T \in \mathbf{R}^m$.
$diag \{A_1, . . . ,A_N\}$ represents
the block diagonal matrix with matrices $A_1, . . . ,A_N$ on its
main diagonal. Denote $col(x_1,....,x_N) $ as the column vector stacked with vectors $x_1,...,x_N$.
$I_n$ denotes the identity matrix in $\mathbf{R}^{n\times n}$.
For a matrix $A=[a_{ij}]$, $a_{ij}$ or $[A]_{ij}$
stands for the matrix entry in the $i$th row and $j$th column of $A$.
Denote $int(\Omega)$ as the interior of $\Omega$ and
$ri(\Omega)$ as the relative interior of $\Omega$.
Denote $\times_{i=1,...,N}\Omega_i$ or $\prod_{i=1}^N \Omega_i$ as the Cartesian product of the sets $\Omega_i,i=1,...,N$. For vectors (matrices) $a$ and $b$, $a\otimes b$ is their Kronecker product, and $a\odot b$ is
their Hadamard product when their dimensions are the same.

Let $\mathfrak{A}:\mathbf{R}^m \rightarrow 2^{\mathbf{R}^m}$ be a set-value operator. Denote ${\rm Id}$ as the identity operator, i.e, ${\rm Id}(x)=x$.
The domain of $\mathfrak{A}$ is $dom\mathfrak{A}= \{x| \mathfrak{A}x \neq \emptyset\}$ where $\emptyset$ stands for the empty set, and the range of $\mathfrak{A}$ is $ran\mathfrak{A}=\{y | \exists x, y\in \mathfrak{A}x\}$. The graph of $\mathfrak{A}$ is $gra\mathfrak{A}=\{(x,u) | u\in \mathfrak{A}x\}$, then the inverse of $\mathfrak{A}$ has its graph as $gra\mathfrak{A}^{-1}=\{(u, x)| (x, u)\in gra \mathfrak{A}\}$.
The zero set of  $\mathfrak{A}$ is $zer\mathfrak{A}=\{x | \mathbf{0} \in \mathfrak{A}x\}$.
The sum of  $\mathfrak{A}$ and $\mathfrak{B}$ is defined as $gra (\mathfrak{A}+\mathfrak{B})=\{(x,y+z)| (x,y)\in gra \mathfrak{A}, (x,z)\in gra \mathfrak{B}\}$. Define the {\it resolvent} of $\mathfrak{A}$ as $R_{\mathfrak{A}}=({\rm Id}+\mathfrak{A})^{-1}$.
The composition of two operators $\mathfrak{A}$ and $\mathfrak{B}$ is defined with its graph $gra\mathfrak{A}\circ \mathfrak{B}= \{(x,u) | \exists v\in \mathfrak{B}(x), u\in \mathfrak{A}(v)\}$.

Operator $\mathfrak{A}$ is called monotone if
$\forall (x,u), \forall(y,v)\in gra\mathfrak{A}$, we have
$\langle x-y, u-v\rangle \geq 0.$
Moreover, it is maximally monotone if $gra\mathfrak{A}$ is not {\it strictly} contained in the graph of any other monotone operator.
A skew-symmetric matrix $A=-A^T$ defines a maximally monotone operator $Ax$ (Example 20.30 of \cite{combettes1}).
Suppose  $\mathfrak{A}$ and $\mathfrak{B}$ are maximally monotone operators and
$0\in int(dom \mathfrak{A}-dom \mathfrak{B})$, then $\mathfrak{A}+\mathfrak{B}$ is also maximally monotone.
For a closed convex set $\Omega$, denote $N_{\Omega}(x)$ as the normal cone operator of  $\Omega$, and
$N_{\Omega}(x)=\{v| \langle v, y-x\rangle\leq 0, \forall y\in \Omega\}  $ and $dom N_{\Omega}=\Omega$. $N_{\Omega}$ is maximally monotone. Define the projection of $x$ onto set $\Omega$ as $P_{\Omega}(x)=\arg\min_{y} || x-y||_2^2$, then $P_{\Omega}(x)=R_{N_{\Omega}}(x)$.

For a single-valued operator $T:\Omega \subset \mathbf{R}^m\rightarrow \mathbf{R}^m$, a point $x\in \Omega$ is a fixed point of $T$ if $Tx=x$, and the set of fixed points of $T$ is denoted as $FixT$.
$T$  is nonexpansive if
it is $1-$Lipschitzian, i.e., $||T(x)-T(y) || \leq ||x-y||, \forall x,y \in \Omega$.
Let $\alpha \in (0,1)$, then $T$ is $\alpha-$averaged  if there exists a nonexpansive operator $T^{'}$ such that  $T=(1-\alpha){\rm Id}+\alpha T^{'}$.
Denote the class of $\alpha-$averaged operators as $\mathcal{A}(\alpha)$.
If $T\in \mathcal{A}(\frac{1}{2})$, then $T$ is called firmly nonexpansive.
An operator $T$ is $\beta-$cocoercive if $\beta T \in \mathcal{A}(\frac{1}{2})$, which satisfies
\begin{equation}\label{equ_cocercive}
\langle x-y,Tx-Ty \rangle\geq \beta ||Tx-Ty ||^2,\forall x,y\in domT.\nonumber
\end{equation}
If operator $\mathfrak{A}$ is maximally monotone, then $T=R_{\mathfrak{A}}=({\rm Id}+\mathfrak{A})^{-1}\in \mathcal{A}(\frac{1}{2})$ and $domR_{\mathfrak{A}}=\mathbf{R}^m$ (Proposition 23.7 of \cite{combettes1}).

\section{Game formulation}\label{sec_game_formulation}

Consider a set of players (agents) $\mathcal{N}=\{1,\cdots, N \}$ involved in the following noncooperative game with shared coupling constraints.
Player $i\in \mathcal{N}$ controls its own decision (strategy or action) $x_i\in \Omega_i \subset \mathbf{R}^{n_i}$, where $\Omega_i $ is its private feasible decision set.  Let  ${x} =col(x_1,\cdots,x_N) \in \mathbf{R}^n$ denote the overall decision profile,  i.e., the stacked vector of all agents' decisions, with $\sum_{i=1}^N n_i=n$.
Let  ${x}_{-i}=col(x_1,\cdots,x_{i-1},x_{i+1},\cdots,x_{N})$ denote the decision profile of all agents except player $i$, then $x=(x_i,x_{-i})$.
 Player $i$ aims to optimize its own objective function within its feasible set,  $f_i(x_i,{x}_{-i}): {\Omega} \rightarrow \mathbf{R}$ where ${\Omega}=\prod_{i=1}^N \Omega_i  \subset \mathbf{R}^{n}$. Note that $f_i(x_i,{x}_{-i})$ is coupled with  other players' decisions ${x}_{-i}$.
Moreover, all players' decisions are coupled together through a {\it globally shared set } $X \subset \mathbf{R}^n$.
Hence, player $i$ has a set-valued map $X_i({x}_{-i}): \mathbf{R}^{n-n_i}\rightarrow 2^{\mathbf{R}^{n_i}}$ that specifies its feasible set  defined as
\begin{equation}\label{feasible_set_map}
X_i({x}_{-i}): =\{x_i\in  \Omega_i| (x_i,{x}_{-i})\in X\}.
\end{equation}
Given ${x}_{-i}$, player $i$'s best-response strategy is
\begin{equation}\label{GM}
\min_{x_i} \;  f_i(x_i,{x}_{-i}),  \; s.t., \;  x_i \in X_i({x}_{-i}).
\end{equation}
A generalized Nash equilibrium (GNE) ${x}^*$ is defined at the intersection of all  players' best-response sets, i.e., $\forall i\in \mathcal{N}$,
\begin{equation}\label{GNE}
x_i^*\in \arg\min_{x_i} f_i(x_i,{x}^*_{-i}),\;  s.t., \;  x_i \in X_i({x}^*_{-i}).
\end{equation}

In this work, we consider the following coupling set,
\begin{equation}\label{coupling_set}
 X := \prod_{i=1}^N \Omega_i \bigcap \{ {x} \in \mathbf{R}^n |  \sum_{i=1}^N A_ix_i = \sum_{i=1}^N b_i\}. 
\end{equation}
where $A_i\in \mathbf{R}^{m\times n_i}$ and $b_i\in \mathbf{R}^m$ as well as $\Omega_i$ are private data of player $i$. 
 Thereby, the shared  set $X$ couples all players' feasible sets, but is {\bf not known} by any agent.
We consider the following  assumption on the game in \eqref{GM}.

{\it \begin{assum}\label{assum1}
$\Omega_i$ is  a closed  convex set with nonempty interiors, and $X$ in \eqref{coupling_set} has nonempty relative interiors.
Given any ${x}_{-i}\in \prod_{j=1,j\neq i}^N \Omega_i$,  the set $X_i({x}_{-i})$ in \eqref{feasible_set_map} has nonempty relative interiors. For player $i$, $f_i(x_i,{x}_{-i})$ is a continuously differentiable convex function with respect to $x_i$ given any fixed ${x}_{-i}$.
\end{assum}}
The {\it pseudo-gradient} of the game  \eqref{GM}, denoted as $F({x})$, is  defined as
\begin{equation}\label{pseudogradient}
F({x})= col(\nabla_{x_1} f_1(x_1,{x}_{-1}),\cdots, \nabla_{x_N} f_N(x_N,{x}_{-N})).
\end{equation}
Then we also impose the following assumption on $F({x})$,
\begin{assum}\label{assum2}
The {\it pseudo-gradient} $F({x})$ of the game  \eqref{GM}
is $\upsilon-$strongly monotone over $\Omega$: $$\langle F({x})-F({y}), {x}-{y}\rangle \geq \upsilon ||{x}-{y}  ||_2^2, \forall {x}, {y}\in \Omega,$$ and
$\chi-$Lipschitz continuous over $\Omega$: $$||  F({x})-F({y})||_2 \leq \chi||{x}-{y}  ||_2,\forall {x}, {y}\in \Omega.$$
\end{assum}

\vskip 3mm
With $F({x})$, we can define the {\it  variational inequality} (VI) problem  as follows,
\begin{equation}\label{vi_original}
Find \; {x}^*, \; s.t.,\;  \langle F({x}^*),{x}-{x}^*\rangle \geq 0, \; \forall {x}\in  X.
\end{equation}
Under Assumption \ref{assum1}, ${x}^*$ is a  GNE of the game in \eqref{GM} if and only if
$\forall i\in \mathcal{N}$ there exists $\lambda_i^* \in \mathbf{R}^m$ such that,
 \begin{equation}
 \begin{array}{l}\label{kkt1}
 \mathbf{0}\in \nabla_{x_i} f_i(x_i^*,{x}^*_{-i})+A_i^T \lambda_i^*+N_{\Omega_i}(x_i^*), \quad \forall i \in
 \mathcal{N},\\
\mathbf{0} = \sum_{i=1}^N A_i x_i^*-\sum_{i=1}^N b_i.
\end{array}
 \end{equation}
 Meanwhile, based on the Lagrangian duality for VI (\cite{pang}), ${x}^*$ is a solution of  VI \eqref{vi_original}  iff  $ \exists \lambda^* \in \mathbf{R}^m $ such that
 \begin{equation}
 \begin{array}{l}\label{kkt2}
  \mathbf{0} \in \nabla_{x_i} f_i(x_i^*,{x}^*_{-i}) + A_i^T \lambda^*+N_{\Omega_i}(x_i^*), \quad \forall i\in \mathcal{N},\\
  \mathbf{0} = \sum_{i=1}^N (A_ix_i^*-b_i).
 \end{array}
 \end{equation}
 By comparing the KKT conditions in \eqref{kkt1} and \eqref{kkt2},  any solution to VI \eqref{vi_original} is a GNE of the game in \eqref{GM}, which is termed as a {\em variational GNE} with all players having the same local  multiplier.

Not every GNE of the  game in \eqref{GM} is a solution to  \eqref{vi_original}.
Since the variational GNE has an economic interpretation of no price discrimination and enjoys a stability  property (refer to \cite{pang} and \cite{yipeng}),  {\it we aim to propose  novel asynchronous distributed algorithms for computing  a variational GNE}.

\begin{rem}
Assumption \ref{assum1} ensures that there exists a unique solution to the VI in \eqref{vi_original}, hence, guarantees the
existence and uniqueness of the variational GNE of the considered game.
\end{rem}

\section{Asynchronous distributed GNE computation with full decision information}\label{sec_algorithm_and_limiting_poinit}

In this section, we propose an asynchronous distributed algorithm with {\it full decision information}, which the players can use to find a solution of VI in \eqref{vi_original} when each player is able to access other players' decisions that its local objective function directly depends on,
hence to compute a variational GNE of game in \eqref{GM}.

 We focus on asynchronous distributed algorithms with delayed information because of following reasons. Firstly,
  player $i$ can only manipulate its local  $f_i(x_i,{x}_{-i})$, $A_i$, $b_i$ and $\Omega_i$ for local computation, since that local data contains its private information.
 Secondly, we assume there is no central node that has bidirectional communications with all players, which is vulnerable  to single point failure.
 Thirdly, it is preferred that the players can compute asynchronously such that they can carry on local iterations without waiting for the slowest agent to finish.
With asynchronous GNE seeking algorithms, there is no global iteration counter, the idle time is removed, and the computational resource of each player is fully exploited.
 The computation can rely on delayed information, and hence,  each agent does not need to wait for  its neighbours in order to carry on its local iteration.

We first introduce the communication graph and related auxiliary variables in Section  \ref{sub_sec_auxiliary_variables}, then give the asynchronous distributed algorithm in Section  \ref{sub_sec_alg}.

\subsection{Communication graph and auxiliary variables}\label{sub_sec_auxiliary_variables}

To facilitate the distributed computation, players are able to communicate with  their neighbours through a {\it connected and undirected} graph $\mathcal{G}=(\mathcal{N},\mathcal{E})$.
The edge set $\mathcal{E} \subset \mathcal{N}\times \mathcal{N} $ contains all the  information interactions.
$(i,j) \in \mathcal{E}$ if agent $i$ and $j$ can share information with each other, and
 $j$, $i$ are  called neighbours.
A path of graph $\mathcal{G}$ is a sequence of distinct agents in
$\mathcal{N}$ such that any consecutive agents in the sequence
are neighbours. Agent $j$ is said to be
connected to agent $i$ if there is a path from $j$ to $i$.
$\mathcal{G}$ is connected since any two agents are
connected.

Obviously, $|\mathcal{N}|=N$, and we denote $|\mathcal{E}|=M$.
The edges are also labeled with $e_l$, $l=1,...,M$.
Without  loss  of  generality, $e_l=(i,j)$ is {\it arbitrarily} ordered as $i\rightarrow j$.
Define $\mathcal{E}^{in}_i$ and $\mathcal{E}^{out}_i$ for  agent $i$ as follows:
$e_l \in  \mathcal{E}^{in}_i$  if agent $i$ is the targeted point of $e_l$; $e_l\in \mathcal{E}^{out}_i $ if agent $i$ is the starting point of $e_l$. Then denote $\mathcal{E}_i=\mathcal{E}^{in}_i \bigcup \mathcal{E}^{out}_i$ as the set of edges adjoint to agent $i$.
Define the incidence matrix of $\mathcal{G}$ as ${V}\in \mathbf{R}^{N\times M}$ with
${V}_{il}=1$ if  $e_l\in \mathcal{E}_i^{in}$, and ${V}_{il}=-1$ if $e_l\in \mathcal{E}_i^{out}$, otherwise ${V}_{il}=0$. We have $\mathbf{1}^T_N V=\mathbf{0}^T_M$,
 and $V^Tx=\mathbf{0}_M$ if and only if $x\in \{\alpha \mathbf{1}_N | \alpha \in\mathbf{ R}\}$ since $\mathcal{G}$ is connected.
 Denote $\mathcal{N}_l=\{i,j\}$ as the pair of agents  connected by edge $e_l=(i,j)$.
We also define the node neighbour set of each agent $i$ as $\mathbb{N}_i=\{j|\mathcal{E}_i \cap \mathcal{E}_j \neq \emptyset \}$.

With the incidence matrix $V$, we can define the {\it Laplacian} of $\mathcal{G}$ as $L=VV^T$ and the {\it edge Laplacian} of $\mathcal{G}$ as $L^e=V^TV$.
It is known that $L_{ij}=-1$ if $j\in \mathbb{N}_i$, and $L_{ii}=-\sum_{j=1,j\neq i}^N L_{ij}$.
Moreover,   $L^e\in \mathbf{R}^{M\times M}$ is  also symmetric, and
$[L^e]_{ll}=2$; $[L^e]_{lp} \neq0 $ for $l \neq p$ if $e_l$ and $e_p$ are connected to the same node $i$, otherwise $[L^e]_{lp}=0$. Moreover,  $[L^e]_{lp}=1$ if $e_l$ and $e_p$  share the
same direction in term of $i$, that is either both point to $i$ or both start from $i$; and $[L^e]_{lp}=-1$ if
$e_l$ and $e_p$ have opposite directions in term of $i$, that is one points to $i$ and the other starts from $i$.
Define an edge neighbour set of $e_l$ as $\mathcal{N}^e_l=\{e_p| [L^e]_{lp}\neq0\}$. Notice that $e_l\in \mathcal{N}^e_l $. Moreover, we know that
$\mathcal{N}^e_l = \bigcup_{j\in \mathcal{N}_l} \mathcal{E}_{j} $. Refer to Chapter 2 of \cite{mugnus} for {\it edge Laplacian} and its application in consensus problems.

We introduce variables in the distributed algorithm.
Firstly, each player controls its local decision $x_i$ and  its {\it local  multiplier} $\lambda_i$. According to KKT \eqref{kkt2}, in the steady state all players should have the same local  multipliers, i.e., $\lambda_i=\lambda^*,\forall i\in \mathcal{N}$.
To facilitate the coordination for the consensus of local multipliers and to ensure the  coupling constraint, we introduce an auxiliary variable $z_l\in \mathbf{R}^m$ related with edge $e_l$ of graph $\mathcal{G}$.
Intuitively speaking, $z_l$ can be regarded as a resource ``flow" on each edge to reach the balance of the equality constraint, hence $z_l$ can be regarded as a dual variable for network flow problems as in  \cite{burger}.
Notice that $\mathcal{G}$ is undirected and  the edges are {\it arbitrarily} ordered, therefore,  any agent from $\mathcal{N}_l$ can  maintain $z_l$.
For clarity,  we consider that $z_l$ will be maintained by agent $i$ if $e_l\in  \mathcal{E}^{out}_i$, but $z_l$  is shared/known by the two players in  $\mathcal{N}_l$.
We assume that each player $i$ has an output buffer where he writes to, which can be read by its neighbouring players $j \in \mathbb{N}_i$,  and an input buffer where he reads from. Since the input buffer also contains $z_l, e_l\in \mathcal{E}^{in}_i$, it can be written by its neighbour $j \in \mathbb{N}_i$, $j\in \mathcal{N}_{l} \backslash i$.

In the algorithm of this section, we assume that player $i$ can read the decision $x_{j}$ from the output buffer of player $j$ if  $f_i(x_{i},{x}_{-i})$ {\it explicitly} depends on $x_j$, and that it has a local oracle that returns $\nabla_{x_i}f_i(x_i,{x}_{-i})$, given ${x}_{-i}$.
Even though in some applications, such as the Nash Cournot game in \cite{sayed} and \cite{yipeng} and the example in Section \ref{sec_numerical_studies},  each player's local objective function only depends on part of the overall decision,  we still call the algorithm as one with {\it full decision information} since each player needs to know all the decisions that its local objective function directly depends on.

\subsection{Asynchronous distributed algorithm with full decision information}\label{sub_sec_alg}

To describe the  \b{A}synchronous  \b{D}istributed \b{A}lgorithm for \b{G}\b{N}\b{E} Computation with Delay\b{s} (ADAGNES), we assume that each agent has a (virtual) Poisson variable (clock), which can model the number of iterations in a given time period. Then the ADAGNES is given as:
\begin{alg}[ADAGNES]\label{alg_d}
\quad \\
\noindent\rule{0.49\textwidth}{0.7mm}

{\bf Initialization:} Player $i$ picks $x_{i,0}\in \Omega_i$, $z_{l,0}\in \mathbf{R}^m$, $ e_l\in \mathcal{E}_i^{out}$  and $\lambda_{i,0}\in \mathbf{R}^m$, and has a local variable $B_{i,0}=A_ix_{i,0}-b_{i,0}$ and a Poisson clock with rata $\varsigma_i$.

\noindent\rule{0.49\textwidth}{0.7mm}
{\bf Iteration at  $k$:}
Suppose player $i_k$'s clock ticks at time $k$, then player $i_k$ is active and updates its local
variables as follows:
\vskip 1mm
\quad {\bf Reading phase}: Get information from its neighbours' output buffer and its input buffer.
Duplicate its input buffer to its output buffer.
Read $x_{j,k-\delta^k_j}$ if $f_{i_k}(x_{i_k},{x}_{-i_k})$ directly depends on $x_j$.
If $ j \in \mathbb{N}_{i_k}$, read $\lambda_{j,k-\delta^k_j}$ from the output buffer of $j$.
If $e_l\in \mathcal{E}^{out}_{i_k}$ and $j \in \mathcal{N}_{l} \backslash i_k$,
then read $B_{j,k-\delta^k_{j}}$ and $z_{q,k-\pi^k_q}$, $e_q \in \mathcal{E}_{j}$ from the output buffer of $j$.

\vskip 1mm

\quad {\bf Computing phase}:

\qquad i): update $\lambda_{i_k}$,
\begin{equation}
\begin{array}{l}\label{equ_delayed_alg_1}
\tilde{\lambda}_{i_k,k}  =  \lambda_{i_k,k-\delta^k_{i_k}}+\sigma( A_{i_k}x_{i_k,k-\delta^k_{i_k}}-b_{i_k} \\
 \quad + \sum_{e_l \in \mathcal{E}_{i_k}}V_{i_kl}z_{l,k-\pi^k_l})\\
\lambda_{i_k,k+1}=\lambda_{i_k,k-\delta^k_{i_k}}+\eta(\tilde{\lambda}_{i_k,k}-\lambda_{i_k,k-\delta^k_i})
\end{array}
\end{equation}

\qquad ii): update $x_{i_k}$,
\begin{equation}
\begin{array}{l}\label{equ_delayed_alg_2}
\tilde{x}_{i_k,k}       =  P_{\Omega_{i_k}}\big[ x_{i_k,k-\delta^k_{i_k}}-\tau(\nabla_{i_k}f_{i_k}(x_{i_k,k-\delta^k_{i_k}},{x}_{-i_k, {k-\delta^k}}) \\
\qquad + A_{i_k}^T(2\tilde{\lambda}_{i_k,k}-\lambda_{i_k,k-\delta^k_{i_k}}))  \big]\\
x_{i_k,k+1}       =  x_{i_k,k-\delta^k_{i_k}}+\eta(\tilde{x}_{i_k,k}-x_{i_k,k-\delta^k_{i_k}})
\end{array}
\end{equation}

\qquad iii): If $e_l\in \mathcal{E}^{out}_{i_k}$ and  $ j_k \in \mathbb{N}_{i_k}$, $j_k \in \mathcal{N}_{l} \backslash i_k$,
then update $z_{l,k}$:
\begin{equation}
\begin{array}{l}\label{equ_delayed_alg_3}
{z}_{l,k+1}        =   z_{l,k-\pi^k_{l}}- \eta\gamma(\lambda_{j_k,k-\delta^k_{j_k}}-\lambda_{i_k,k-\delta^k_{i_k}})\\
\qquad-2\eta\sigma\gamma\big[(B_{j_k,k-\delta^j_{j_k}}-B_{i_k,k-\delta^k_{i_k}})+\sum_{q\in \mathcal{N}^e_{l}}L^e_{lq}z_{q,k-\pi^k_{q}}\big]
\end{array}
\end{equation}

\vskip 1mm
\quad {\bf Writing phase}:
 Write $x_{i_k,k+1}$, $\lambda_{i_k,k+1}$, $z_{l,k+1},e_l
 \in \mathcal{E}^{out}_{i_k}$  and $B_{i_k,k+1}=A_{i_k}x_{i_k,k+1}-b_{i_k}$ to its output buffer, and write $z_{l,k+1}$, $e_l\in \mathcal{E}^{out}_{i_k}$ to the  input buffer of player $j_k \in \mathcal{N}_{l} \backslash i_k$.

\vskip 3mm

All other players keep their variables (output buffers) unchanged, and only increase $k$ to $k+1$.

\noindent\rule{0.49\textwidth}{0.7mm}
In the foregoing,  $\sigma$, $\gamma$, $\tau$ and $\eta$ are real positive step-sizes, and
$x_{i,k}, \lambda_{i,k}$ and $z_{l,k}$ are $x_i,\lambda_i$ and $z_l$, respectively, at iteration $k$.
$\delta^k_i$ and $\pi^k_l$ are the delays of $x_i$ (as well as $\lambda_i$ and $B_i$) and $z_l$, respectively, to characterize
the delayed information used by the active player $i$ at time $k$.
${x}_{-i, {k-\delta^k}}$ is the vector composed with the delayed decision ${x}_{j, {k-\delta^k_j}}$ which  directly influences  $f_i(x_i,{x}_{-i})$ of the active player $i$ at time $k$.
\end{alg}

\vskip 2mm
\begin{rem}
Player $i$ can get the $l-$th row of  $L^e$ if $e_l\in \mathcal{E}_i^{out} $ with one step communication before the  algorithm starts.
With the duplication step in reading phase, the input buffer of agent $i$ contains $z_l, e_l\in \mathcal{E}^{in}_i$, while its output buffer contains
$x_{i_k,k+1}$, $\lambda_{i_k,k+1}$, $z_{l,k+1},e_l
 \in \mathcal{E}^{out}_{i_k}$ as well as $z_l, e_l\in \mathcal{E}^{in}_i$.
Each player $i$ can get
$z_{q,k-\pi^k_q}, e_q\in \mathcal{N}^e_{l}$, $e_l\in \mathcal{E}_i^{out}$ by just inquiring $z_{p,k-\pi^k_p}$, $e_p \in \mathcal{E}_{j}$ from its neighbour $j=\mathcal{N}_l\backslash i$. Hence the algorithm is  distributed without any center.
The assumption that each player could read all  decisions that influence its $f_i(x_{i},{x}_{-i})$ might be restrictive in some applications, such as in aggregative games.
This will be addressed by introducing local estimations of the overall decision profile and incorporating an additional local consensus dynamics in Section \ref{sec_pdi}.
\end{rem}

\begin{rem}
Since agent $i$ updates $x_i$, $\lambda_i$, $B_i$ and $z_l, e_l \in \mathcal{E}^{out}_i $ simultaneously, those variables share the same delay at $k$, i.e., $\delta^k_i=\pi^k_l, e_l\in \mathcal{E}^{out}_i$. However, the variables $z_l,e_l\in \mathcal{E}^{out}_i$ and $z_l,e_l\in \mathcal{E}^{in}_i$ may have different delays, since they are maintained by different agents. Thereby, both $\delta^k_i$ and $\pi^k_l$ are introduced.
Moreover, if player $i_k$ is active at time $k$, we have $x_{i_k,k-\delta^k_{i_k}}=x_{i_k,k-\delta^k_{i_k}+1}=,\cdots,=x_{i_k,k}$, while similar relations also hold for $\lambda_{i}$ and $z_{l}, e_l\in \mathcal{E}_i^{out}$. To explicitly show how the iterate depends on the delayed information, we will use the latest index when the coordinate is updated.
\end{rem}

\begin{rem}
The delays arise from various sources. It is  possible that one agent has finished the writing phase (index $k$ will increase 1) before another agent finishes its computation, in which case the index of the information used by the agents in computation phase gets delayed by one. Overall, the delays account for the time due to reading buffers, computation with local data and writing buffers. The algorithm is  asynchronous since the agent needn't wait for its neighbours before carrying on its local iteration.
\end{rem}

\begin{rem}
Each round of iteration is composed of three phases: reading information from neighbour  players' buffers, computation with local data and the delayed information, writing the updated information to the  buffers.
The iteration index $k$ increases only if an agent  has finished its writing phase.
We assume that a buffer is ``locked" if an agent is reading from or writing to that buffer.
We adopt a Poisson clock model since usually the total time consumed for reading, computing and writing phases can be modeled as a random variable with an exponential distribution.
\end{rem}

\section{Algorithm development and analysis}\label{sec_convergence_analysis}

In this section, we first show how the algorithm without asynchrony and delays is developed based on an operator theoretic approach,
and then give the convergence analysis for the asynchronous one with full decision information based on fixed-point iterations in \cite{yinwotao_arock}.

\subsection{Algorithm development without asynchrony and delays}\label{sub_sec_expand}

Let us consider the {\it \b{S}\b{y}nchronous \b{D}istributed GE\b{N} Computation Algorithm without D\b{e}la\b{y}s}, which is denoted as SYDNEY. In this case, the update of agent $i$ at time $k$ is given as:
\begin{alg}[SYDNEY]\label{alg_d_without_delay}
\begin{equation}
\begin{array}{l}\label{equ_syn_no_delay}
\tilde{\lambda}_{i,k}  =  \lambda_{i,k}+\sigma( A_{i}x_{i,k}  + \sum_{e_l \in \mathcal{E}_{i}}V_{il}z_{l,k}-b_{i})\\
\lambda_{i,k+1}=\lambda_{i,k}+\eta(\tilde{\lambda}_{i,k}-\lambda_{i,k})\\
\tilde{x}_{i,k}       =  P_{\Omega_{i}}\big[ x_{i,k}-\tau(\nabla_{i}f_{i}(x_{i,k},{x}_{-i, k})
+ A_i^T(2\tilde{\lambda}_{i,k}-\lambda_{i,k}))  \big]\\
x_{i,k+1}       =  x_{i,k}+\eta(\tilde{x}_{i,k}-x_{i,k})\\
{z}_{l,k+1}        =   z_{l,k}- \eta\gamma(\lambda_{j,k}-\lambda_{i,k})-2\eta\sigma\gamma(B_{j,k}-B_{i,k})\\
\quad -2 \eta\sigma\gamma\sum_{q\in \mathcal{N}^e_{l}}L^e_{lq}z_{q,k} \nonumber
\end{array}
\end{equation}
The notations are the same as  in ADAGNES.
\end{alg}

{Here are the compact notations.
 Denote ${x}_{k}=col(x_{1,k},\cdots,x_{N,k})$ as all players' decisions, and
 $\tilde{{x}}_k=col(\tilde{x}_{1,k},\cdots,\tilde{x}_{N,k})$.
 Denote ${{\Lambda}}_k=col(\lambda_{1,k},\cdots,\lambda_{N,k})$, and $\tilde{{\Lambda}}_k$ are defined similarly.
Denote $\mathbf{z}_k=col(z_{l,k})_{l\in \mathcal{E}}$, and we also introduce an intermediate variable
$\tilde{\mathbf{z}}_k=col(\tilde{z}_{l,k})_{l\in \mathcal{E}}$.
Denote $\mathbf{b}=col(b_1,\cdots,b_N)$, $\mathbf{A}= \diag\{A_1,\cdots,A_N\}$, $\mathbf{V}=V\otimes I_m$  and
$\mathbf{L}^e=L^e\otimes I_m$.}
Then SYDNEY  can be written in a compact form as,
\begin{equation}
\begin{array}{l}\label{equ_syn_compact}
\tilde{\Lambda}_{k} = {\Lambda}_{k}+\sigma(\mathbf{A} \mathbf{x}_{k} +  \mathbf{V} \mathbf{z}_k-\mathbf{b}), \\
\tilde{\mathbf{z}}_{k}    = \mathbf{z}_{k}-\gamma \mathbf{V}^T {\Lambda}_{k} - 2\sigma\gamma\mathbf{V}^T (\mathbf{A}  \mathbf{x}_{k}-\mathbf{b})
-2 \sigma\gamma \mathbf{L}^e\mathbf{z}_k, \\
\tilde{{x}}_{k}    =P_{\Omega}\big[   {x}_k-\tau(F({x})+\mathbf{A}^T(2\tilde{\Lambda}_{k}-\Lambda_{k}))\big],\\
{\Lambda}_{k+1}  =   \Lambda_k  +  \eta(\tilde{\Lambda}_{k}-\Lambda_k ), \quad
{\mathbf{z}}_{k+1}     =   \mathbf{z}_{k}   +  \eta (\tilde{\mathbf{z}}_{k}- \mathbf{z}_{k}),\\
{{x}}_{k+1}     =   {x}_k     +  \eta (\tilde{{x}}_{k}- {x}_k  ).
\end{array}
\end{equation}

Next we show that \eqref{equ_syn_compact} is derived by a preconditioned forward-backward splitting algorithm for finding zeros of a sum of monotone operators. Denote
$\varpi=col(\Lambda,\mathbf{z},{x})$, and then
define operators $\mathfrak{A}$ and $\mathfrak{B}$ as follows,
\begin{equation} \label{operator_A_B}
{\mathfrak{A}}: \varpi \mapsto
\left(
  \begin{array}{c}
    -\mathbf{A}\mathbf{x} -\mathbf{V}\mathbf{z}    \\
    \mathbf{V}^T\Lambda   \\
    \mathbf{A}^T\Lambda+  N_{\Omega}(x) \\
  \end{array}
\right),{\mathfrak{B}}: \varpi \mapsto
\left(
  \begin{array}{c}
 \mathbf{b}    \\
   \mathbf{0 } \\
   F({x}) \\
  \end{array}
\right)
\end{equation}

\begin{thm}\label{thm_limiting_point_ine_1}
Suppose that Assumption \ref{assum1}-\ref{assum2} hold for game \eqref{GM}.
Then, any zero $col(\Lambda^*,\mathbf{z}^*,{x}^*)$ of ${\mathfrak{A}}+\mathfrak{B}$ in \eqref{operator_A_B} has the ${x}^*$ component as a variational GNE of \eqref{GM}, and $\Lambda^*=\mathbf{1}_N\otimes \lambda^*$. ${x}^*$ and $\lambda^*$ satisfy KKT condition \eqref{kkt2}.
\end{thm}

Theorem \ref{thm_limiting_point_ine_1} relates the zero of ${\mathfrak{A}}+\mathfrak{B}$ to the variational GNE of the game. In fact,
$\mathbf{0}\in ({\mathfrak{A}}+\mathfrak{B})\varpi^*$ can be regarded as a decomposition of KKT \eqref{kkt2} by introducing local multipliers $\lambda_i$ and edge variables $z_l,e_l\in \mathcal{E}$.
The proof of Theorem \ref{thm_limiting_point_ine_1} utilizes both the properties $ker{V^T}=span\{\mathbf{1}_N\}$ when  $\mathcal{G}$ is connected and $\mathbf{1}^T V=\mathbf{0}$. Since related and similar analysis can be found in Theorem 4.5 of \cite{yipeng} and Theorem 1 of \cite{yipeng2},  its proof is omitted here.

We show next that  SYDNEY \eqref{equ_syn_compact} can be further written as
\begin{equation}\label{fixed point iteration}
\varpi_{k+1}=\varpi_k + \eta(T\varpi_k - \varpi_k).
\end{equation}
where  $\varpi_k=col(\Lambda_k,\mathbf{z}_k,{x}_k)$, and  the matrix $\Phi$ and operator $T$ are defined as follows,
\begin{equation}\label{precon_phi}
\Phi=\left(
       \begin{array}{ccc}
         \sigma^{-1}I & \mathbf{V}       & \mathbf{A} \\
         \mathbf{V}^T    & \gamma^{-1}I  & \mathbf{0} \\
         \mathbf{A}^T    & \mathbf{0}    & \tau^{-1}I \\
       \end{array}
     \right),
\end{equation}
\begin{equation}\label{operator_T}
T= ({\rm Id}-\Phi^{-1}\mathfrak{A} )^{-1}({\rm Id}-\Phi^{-1}\mathfrak{B}).
\end{equation}
Suppose $\sigma$, $\gamma$ and $\tau$ are chosen such that $\Phi$ in \eqref{precon_phi} is positive definite.
Then, note that we have $FixT=zer(\mathfrak{A}+\mathfrak{B})$, and Theorem \ref{thm_limiting_point_ine_1} applies. In fact, $\varpi=T\varpi \Leftrightarrow \Phi\varpi-\mathfrak{B}\varpi\in \mathfrak{A}{\varpi} +\Phi\varpi \Leftrightarrow \mathbf{0}\in (\mathfrak{A}+\mathfrak{B})\varpi$.

To show that SYDNEY \eqref{equ_syn_compact} can be written as \eqref{fixed point iteration},
first denote $\tilde{\varpi}_k=col(\tilde{\Lambda}_k, \tilde{\mathbf{z}}_k,\tilde{{x}}_k)$, and write \eqref{fixed point iteration} as
$\tilde{\varpi}_k=T\varpi_k $ and
$\varpi_{k+1}=\varpi_{k}+\eta(\tilde{\varpi}_k-\varpi_k)$. Now, $\varpi_{k+1}=\varpi_{k}+\eta(\tilde{\varpi}_k-\varpi_k)$ expands to give the  last three equations in  \eqref{equ_syn_compact}.
Moreover, using \eqref{operator_T},  $\tilde{\varpi}_k=T\varpi_k $ is equivalent to
\begin{equation}\label{equ_forward_backward_splitting}
-\mathfrak{B}\varpi_{k}\in \mathfrak{A}\tilde{\varpi}_{k} +\Phi(\tilde{\varpi}_{k}-\varpi_k).
\end{equation}
Using the definition of $\mathfrak{A}$ and $\mathfrak{B}$ in \eqref{operator_A_B} and
$\tilde{\varpi}_k=col(\tilde{\Lambda}_k, \tilde{\mathbf{z}}_k,\tilde{{x}}_k)$, we show next that by expanding \eqref{equ_forward_backward_splitting}  we can can obtain the first three equations in \eqref{equ_syn_compact}  as follows.

The update of ${\Lambda}$ component in  \eqref{equ_forward_backward_splitting} gives
$$
-\mathbf{b}= -\mathbf{A} \tilde{{x}}_k-\mathbf{V}\tilde{\mathbf{z}}_k + \sigma^{-1}({\tilde{\Lambda}}_k-{\Lambda}_k)+\mathbf{A} (\tilde{{x}}_k-{x}_k)+\mathbf{V}(\tilde{\mathbf{z}}_k-\mathbf{z}_k), \nonumber
$$
which yields the first line of \eqref{equ_syn_compact}.

The update of $\mathbf{z}$ component in \eqref{equ_forward_backward_splitting} reads as
$$
\mathbf{0}= \mathbf{V}^T{\tilde{\Lambda}} + \mathbf{V}^T({\tilde{\lambda}}_k-\bar{\lambda}_k )+ \gamma^{-1}({\tilde{\mathbf{z}}}_k-{\mathbf{z}}_k),  \nonumber
$$
hence,
$$
{\tilde{\mathbf{z}}}_k={\mathbf{z}}_k-\gamma\mathbf{V}^T(2{\tilde{\Lambda}}_k-{\Lambda}_k). \nonumber
$$
However, from the first line of \eqref{equ_syn_compact} we have $2{\tilde{\Lambda}}_k-{\Lambda}_k= {\Lambda}_{k}+2\sigma(\mathbf{A} {x}_{k} +  \mathbf{V} \mathbf{z}_k-\mathbf{b})$.
Then,
$$
{\tilde{\mathbf{z}}}_k={\mathbf{z}}_k-\gamma\mathbf{V}^T({\Lambda}_{k}+2\sigma(\mathbf{A} {x}_{k} +  \mathbf{V} \mathbf{z}_k-\mathbf{b})), \nonumber
$$
which yields the second line of \eqref{equ_syn_compact} by  using $L^e=V^TV$.

The update of ${x}$ component  in \eqref{equ_forward_backward_splitting} is
$$
-F({x}_k) \in   \mathbf{A}^T{\tilde{\Lambda}}_k+  N_{{\Omega}}(\tilde{{x}}_k) + \mathbf{A}^T({\tilde{\Lambda}}_k-{\Lambda}_k)
+\tau^{-1}(\tilde{{x}}_k-{x}_k). \nonumber
$$
Using $P_{\Omega}(x)=R_{N_{\Omega}}(x)=({\rm Id} + N_{\Omega})^{-1}x$, we  get
the third line of \eqref{equ_syn_compact}. Thus SYDNEY \eqref{equ_syn_compact} can be written as \eqref{fixed point iteration}.

Moreover, \eqref{fixed point iteration} is an iteration for  computing fixed points of $T$, and any fixed point of
$T$ is also a zero of $\mathfrak{A}+\mathfrak{B}$, which satisfies Theorem  \ref{thm_limiting_point_ine_1}. Under Assumption \ref{assum1}, it can be shown that $T$ is an averaged operator under a proper chosen norm with the following Lemma \ref{lem_operator_A_B} and Lemma \ref{lem_preconditioned_operators}.

\begin{lem}\label{lem_operator_A_B}
Suppose Assumption \ref{assum1}-\ref{assum2} hold. Operator $\mathfrak{A}$ in \eqref{operator_A_B} is maximally monotone. Operator $\mathfrak{B}$ in \eqref{operator_A_B} is $\frac{\upsilon}{\chi^2}$-cocoercive;
\end{lem}
{\bf Proof:}
$\mathfrak{A}$ is the sum of a skew-symmetric linear operator and $\mathbf{0} \times \mathbf{0} \times N_{{\Omega}}(x)$, which are both maximally monotone by Example 20.30, Theorem 20.40 and Example 16.12 of \cite{combettes1}. Hence, $\mathfrak{A}$ is maximally monotone.
For operator $\mathfrak{B}$, by the strongly monotonicity  and Lipschitz continuity of $F({x})$ we have
$ \langle \mathfrak{B}(\varpi_1)-\mathfrak{B}(\varpi_2),\varpi_1-\varpi_2\rangle=\langle F({x}_1)-F({x}_2), {x}_1-{x}_2\rangle \geq \upsilon||{x}_1-{x}_2 ||^2 \geq
\frac{\upsilon}{\chi^2}||F({x}_1)-F({x}_2) ||^2=\frac{\upsilon}{\chi^2}||\mathfrak{B}(\varpi_1)-\mathfrak{B}(\varpi_2) ||^2$, which implies that
$\mathfrak{B}$ in \eqref{operator_A_B} is $\frac{\upsilon}{\chi^2}$-cocoercive by Definition 4.4 and Remark 4.24 in \cite{combettes1}.
\hfill $\Box$

\begin{lem}\label{lem_preconditioned_operators}
Suppose Assumption\ref{assum1}-\ref{assum2} hold.
Take $ \delta > \frac{\chi^2}{2\upsilon}$, and take the step-sizes $\tau, \gamma, \sigma$ in Algorithm \ref{alg_d}  such that $\Phi-\delta I$ is positive semi-definite.
Then under the $\Phi-$induced norm $||\cdot||_{\Phi}$, we have,
(i): $\Phi^{-1}\bar{\mathfrak{A}}$ is  maximally monotone, and $T_1=R_{\Phi^{-1}\bar{\mathfrak{A}}}=({\rm Id}+\Phi^{-1}\bar{\mathfrak{A}})^{-1}\in \mathcal{A}(\frac{1}{2})$.
(ii): $\Phi^{-1}\bar{\mathfrak{B}}$ is $ \frac{\delta\upsilon}{\chi^2}-$cocoercive, and  $T_2={\rm Id}-\Phi^{-1}\bar{\mathfrak{B}}\in \mathcal{A}(\frac{\chi^2}{2\delta{\upsilon}{}})$.
(iii): $T$ in \eqref{operator_T} is  an averaged operator, and $T \in \mathcal{A}(\frac{2\delta \upsilon}{4\delta \upsilon-\chi^2}) $.
\end{lem}
{\bf Proof:} The proof of (i) and (ii) follow from Lemma 5.6 of \cite{yipeng}, and  (iii) follows from Proposition 2.4 of \cite{combettes3}.
\hfill $\Box$

Hence, the convergence of \eqref{fixed point iteration} follows from the Krasnosel'ski\u{i}-Mann (K-M) algorithm
by choosing $0<\eta < \frac{4\delta \upsilon-\chi^2}{2\delta \upsilon}$ according to Theorem 5.14 of \cite{combettes1}.

\begin{rem}
The auxiliary edge variables are also used in \cite{yipeng2} for GNE seeking of monotone games  based on  {\it preconditioned proximal algorithms}.
Moreover, \cite{yipeng2} adopts a preconditioning matrix similar to $\Phi$ in \eqref{precon_phi} to give a {\it two-time-scale or double-layer} algorithm.
However, the algorithm in \cite{yipeng2} is synchronous and involves solving subproblems at each iteration due to the  relaxed monotonicity assumption. On the other hand,
Algorithm \ref{alg_d_without_delay}  is motivated by {\it preconditioned forward-backward splitting algorithms}. Similar ideas have been used in \cite{yipeng} to
give a synchronous GNE seeking algorithm under the same strong monotonicity assumption.
However, the algorithm in \cite{yipeng}  works in a Gauss-Seidel updating order, hence is not implementable in an asynchronous manner.
The key observation of this work is that the combination of {\it edge variables and edge Laplacian  with operator splitting methods
gives a distributed GNE seeking algorithm that can work asynchronously with delayed information}.
\end{rem}

\subsection{Convergence analysis for the asynchronous ADAGNES with full decision information}\label{sec_sub_har}

In the previous subsection, we have shown that the SYDNEY, the synchronous version of ADAGNES, can be written compactly as a fixed-point iteration \eqref{fixed point iteration}. In this subsection, we will give the convergence analysis of the asynchronous ADAGNES,  Algorithm \ref{alg_d} by treating it as  a randomized block-coordinate fixed-point iteration with delayed information, and using results in \cite{yinwotao_arock}.

Denote $\varpi^i=col( x_i, \lambda_i, \{z_l\}_{e_l\in \mathcal{E}_i^{out}}  )$ as the coordinates that should be updated by player $i$.
Define a group of vectors $\Xi_1,\cdots, \Xi_N \in R^{n+Nm+Mm}$, such that for $j=1,\cdots,n+Nm+Mm$,
$[\Xi_i]_j=1$ if the $j-$th coordinate of $\varpi$ is also a coordinate of $\varpi^i$,
and $[\Xi_i]_j=0$, otherwise.
Denote $\xi$ as a random variable that takes values in $\{\Xi_1,\cdots, \Xi_N \}$, and $\mathbb{P}(\xi=\Xi_i)=\varsigma_i/ \sum_{j=1}^N \varsigma_j$, where $\varsigma_i,i=1,\cdots, N$ are the rates of the Poisson clock in Algorithm \ref{alg_d}.
Then let $\xi_k,k=1,...$ to be a sequence of i.i.d. random variables
 that are all equal to $\xi$.
A randomized block-coordinate fixed-point iteration for \eqref{fixed point iteration}, or equivalently for SYDNEY, is
\begin{equation}\label{equ_randomized}
\varpi_{k+1}=\varpi_{k} + \eta \xi_k \odot (T\varpi_k-\varpi_k).
\end{equation}
In \eqref{equ_randomized},  only one agent $i$ is activated at each iteration, with probability $\varsigma_i/ \sum_{j=1}^N \varsigma_j$, which updates using SYDNEY.
 The convergence of \eqref{equ_randomized} with nonexpansive $T$ is fully investigated in \cite{combettes2}.
However, \eqref{equ_randomized} does not account for any delays. Rather it assumes that the operator $T$ is applied to $\varpi_k$, which is the whole vector of all agents' iterates at the current time.
This means that the communication and computation at each iteration
should be performed ``almost instantly", or that players should wait until the current iteration finishes before the next activation, such that each iteration is computed with the most recent information.

To remove this impractical assumption and to reduce the idle time of the players, we allow for delayed information when evaluating $T$ in \eqref{equ_randomized}, which yields
\begin{equation}\label{equ_randomized_delayed}
\varpi_{k+1}=\varpi_{k} + \eta \xi_k \odot (T\hat{\varpi}_k - \varpi_k),
\end{equation}
where $\hat{\varpi}_k$ is the delayed information at time $k$. Next, we show the convergence of ADAGNES by showing that it can be written as  \eqref{equ_randomized_delayed} for appropriately defined $\hat{\varpi}_k$, and by resorting to the theoretical result in \cite{yinwotao_arock}, on asynchronous coordinate updates of  K-M iteration for  nonexpansive operators.

\begin{assum}\label{assum3}
The maximal delay at each step is bounded by $\Psi$, i.e., $\sup_{k} \max_{i\in \mathcal{N}, l\in \mathcal{E}} \{ \max\{\delta_i^k,\pi^k_l\}\} \leq \Psi$.
\end{assum}

\begin{thm}\label{thm_convergence_asy}
Suppose Assumption \ref{assum1}, \ref{assum2} and \ref{assum3} hold.
Take $ \delta > \frac{\chi^2}{2\upsilon}$, and take the step-sizes $\tau, \gamma, \sigma$ in Algorithm \ref{alg_d}  such that $\Phi-\delta I$ is positive semi-definite. Denote $p_{min}=\min_{i\in \mathcal{N}}\{\frac{\varsigma_i}{\sum_{j=1}^N \varsigma_j}\}$, and choose $0<\eta \leq  \frac{c N p_{\min}}{2\Psi\sqrt{p_{\min}}+1} \frac{4\delta \upsilon -\chi^2}{2\delta \upsilon}$ for any $c\in (0,1)$.
Then with Algorithm \ref{alg_d}, the players' decisions ${x}_{k}$ converge to the variational GEN ${x}^*$ of game in \eqref{GM}, and the players' local multipliers $\lambda_{i,k},i\in \mathcal{N}$ reach consensus and  converge to the same $\lambda^*$ with probability $1$.
\end{thm}

{\bf Proof:}
{\it Step 1:} We first show that ADAGNES,  Algorithm \ref{alg_d} can be written as the iteration \eqref{equ_randomized_delayed}  with properly defined delayed information $\hat{\varpi}_k$.
Suppose that player $i_k$ is activated at time $k$, then $\hat{\varpi}_k$ could be constructed as follows.
For player $i_k$ and player $j \in \mathbb{N}_{i_k}$, the corresponding
$x_{i_k}$, $x_j$ and $\lambda_{i_k}$, $\lambda_j$ components of $\hat{\varpi}_k$ are $x_{i_k,k-\delta^k_{i_k}}$, $x_{j,k-\delta^k_j}$ and $\lambda_{i_k, k-\delta^k_{i_k}}$, $\lambda_{j, k-\delta^k_j}$. For a player $j$, whose decision directly influences $f_{i_k}(x_{i_k}, {x}_{-i_k})$, the corresponding $x_j$ component of $\hat{\varpi}_k$ is $x_{j,k-\delta^k_j}$, while for any other player $q$, the corresponding $x_q$ and $\lambda_q$ components of $\hat{\varpi}_k$ can be any element in $\{x_{q,k-\delta^k_q},\cdots x_{q,k}\}$ and $\{\lambda_{j,k-\delta^k_q},\cdots, \lambda_{q,k}\}$ with a bounded $\delta^k_q\leq \Psi$. (In fact, those elements will not directly influence the update of player $i_k$).
For edge $ e_l \in \bigcup_{p\in \mathcal{E}_{i_k}} \mathcal{N}^e_p  $, the corresponding $z_l$
components of $\hat{\varpi}_k$ are $z_{l,k-\pi^k_l}$, while the $z_q$ components of $\hat{\varpi}_k$ corresponding to the other edge $q$    could be any element in $\{z_{q,k-\pi^k_q},\cdots z_{q,k}\}$ with a bounded $\pi^k_q\leq \Psi$.
With the above definition of $\hat{\varpi}_k$, one can expand \eqref{equ_randomized_delayed} as in Section \ref{sub_sec_expand} and obtain the updates of $i_k$ exactly in  equations \eqref{equ_delayed_alg_1}-\eqref{equ_delayed_alg_3}. In fact, by comparing the synchronous ones in Algorithm \ref{alg_d_without_delay} with \eqref{equ_delayed_alg_1}-\eqref{equ_delayed_alg_3}, we only need to replace the $\varpi_k$ in \eqref{equ_randomized} with $\hat{\varpi}_k$ and utilize the fact that $x_{i_k,k-\delta^k_{i_k}}=x_{i_k,k-\delta^k_{i_k}+1}=,\cdots,=x_{i_k,k}$ and the similar relations for $\lambda_{i_k}$ and $z_{l}, e_l\in \mathcal{E}_{i_k}^{out}$.
With the random variable $\xi_k$, the components of $\varpi_k$ except $\varpi^{i_k}$ are kept the same as the previous time. Hence, Algorithm \ref{alg_d} can be written as \eqref{equ_randomized_delayed}.

\vskip 3mm
{\it Step 2:}
By Lemma \ref{lem_preconditioned_operators} and the definition of averaged operators,  we can find a nonexpansive operator $\mathfrak{R}$ such that
$T=(1-\frac{2\delta \upsilon}{4\delta \upsilon-\chi^2}){\rm Id}+ \frac{2\delta \upsilon}{4\delta \upsilon-\chi^2} \mathfrak{R}$. Obviously, $Fix T=Fix \mathfrak{R}$. Therefore,
 the compact form of Algorithm \ref{alg_d}, i.e. \eqref{equ_randomized_delayed}, can be written as $$ \varpi_{k+1}=\varpi_{k} + \eta\xi_k\odot [ \hat{\varpi}_k- \varpi_k + \frac{2\delta \upsilon}{4\delta\upsilon-\chi^2}(\mathfrak{R}\hat{\varpi}_k-\hat{\varpi}_k)].$$
 With $x_{i_k,k-\delta^k_{i_k}}=x_{i_k,k-\delta^k_{i_k}+1}=,\cdots,=x_{i_k,k}$ and the similar relations for $\lambda_{i_k}$ and $z_{l}, e_l\in \mathcal{E}_{i_k}^{out}$, we have
 $\xi_k\odot ( \hat{\varpi}_k- \varpi_k)=0$, i.e., $\xi_k\odot\hat{\varpi}_k=\xi_k\odot \varpi_k$. Hence, Algorithm \ref{alg_d} can be equivalently written as:
\begin{equation}\label{equ_nonexpansive_algorithm}
\varpi_{k+1}=\varpi_{k} +  \frac{2\eta\delta \upsilon}{4\delta \upsilon-\chi^2} \xi_k \odot (\mathfrak{R}\hat{\varpi}_k-\varpi_k).
\end{equation}
According to the Poisson clock assumption, the probability that player $i$ is actived at time $k$ is $\mathbb{P}(\xi_k=\Xi_i)=\frac{\varsigma_i}{\sum_{j=1}^N \varsigma_j}$. Moreover, $\xi_k, k=1,\cdots$ are i.i.d random variables.
Applying { Lemma 13 and Theorem  14 of \cite{yinwotao_arock}}, it follows that $\varpi_k$ generated by  \eqref{equ_nonexpansive_algorithm}  is bounded and, if the step-size $\eta$ is chosen such that
$\frac{2\eta\delta \upsilon}{4\delta \upsilon-\chi^2} \in (0, \frac{cNp_{min}}{\Psi \sqrt{p_{min}}+1}]$ for $c\in (0,1)$, $\varpi_k$ converges to a random variable that takes value in the fixed points of $\mathfrak{R}$. Recall that $Fix \mathfrak{R} = Fix T$ and since  $Fix T=zer(\mathfrak{A}+\mathfrak{B})$, it follows that  $\varpi_k$ generated by Algorithm \ref{alg_d} converges to a random variable that takes values  in  $zer(\mathfrak{A}+\mathfrak{B})$ almost surely. Since the considered game has a unique variational GNE and, by Theorem \ref{thm_limiting_point_ine_1}, any zero $\varpi^*$ of $zer(\mathfrak{A}+\mathfrak{B})$ has its $x^*$ as a variational GNE  and ${\Lambda}^*=\mathbf{1}_N\otimes \lambda^*$,  the conclusion follows.
\hfill $\Box$

\section{Asynchronous algorithm with partial  decision information}\label{sec_pdi}

The asynchronous Algorithm \ref{alg_d} is distributed in the sense that each player only needs to know its local data $f_i,A_i,b_i,\Omega_i$ and to communicate with neighbors.
However, it requires each player $i$ to be able to access all other players' decisions that its $f_i$ directly depends on.
This is  restrictive in  applications. For example, if $f_i(x_i,x_{-i})$ depends on all other players' decisions, like the aggregative games in \cite{shanbhag1,grammatico_2,liangshu}, Algorithm \ref{alg_d} requires  player $i$  to be able to
know all players' (delayed) decisions. Therefore,  Algorithm \ref{alg_d} is called an algorithm  with {\it full decision information}.
Motivated by \cite{pavel3,pavel5,tatarenko,shanbhag1,liangshu,hu}, we show that when the objective functions $f_i$ satisfy additional assumptions,
the proposed framework can lead to an asynchronous GNE computation algorithm even when each player is not able to fully know all other players' decisions that its $f_i$ directly depends on.

\subsection{Asynchronous GNE seeking algorithm with partial decision information}\label{sub_sec_pdi_algorithm}

The key idea is to introduce a local estimation of the overall decision profile for each player, motivated by \cite{hu1,hu2,pavel3,pavel4,wanglong} etc..
We assume that each player $i$ has a local estimation of the GNE ${x}^*$ as
$${x}^{(i)}=[  x^{(i)}_1, \cdots, x^{(i)}_{i-1}, x_i, x^{(i)}_{i+1}, \cdots,x^{(i)}_N] \in \mathbf{R}^{n}.$$
Here $x_i$ is just the local decision of player $i$, and
${x}^{(i)}_{-i}=[  x^{(i)}_1, \cdots, x^{(i)}_{i-1}, x^{(i)}_{i+1}, \cdots,x^{(i)}_N  ]$ is its estimation of all other player's decision. More specific, $x^{(i)}_{j},j\neq i$ is player $i$'s estimation of the decision of player $j$.
We assume that each player $i$ is able to inquire only the local estimations ${x}^{(j)}$ of its neighbors $j\in \mathbb{N}_i$ through the communication graph $\mathcal{G}$ defined in subsection \ref{sub_sec_auxiliary_variables}.  In this case, each player $i$ is not able to evaluate $\nabla_i f_i$ with the (delayed) true decision $(x_i,{x}_{-i})$ but only with its local estimation $(x_i,{x}^{(i)}_{-i})$. Hence, ${x}^{(i)}_{-i}$ should also be updated accordingly to track the true ${x}_{-i}$.
All other local variables and graph notations,  like local multiplier $\lambda_i$, edge variable $z_l$, incidence matrix $V$, Laplacian $L$ and edge Laplacian $L^e$, are defined the same as subsection \ref{sub_sec_auxiliary_variables}.
Then the proposed asynchronous algorithm with partial decision information is given as follows.

\begin{alg}[ADAGNES-PDI]\label{alg_d_partial information}
\quad \\
\noindent\rule{0.49\textwidth}{0.7mm}

{\bf Initialization:} Player $i$ picks $x_{i,0}\in \Omega_i$, a local estimation of other players' decision ${x}^{(i)}_{-i,0}$,  $z_{l,0}\in \mathbf{R}^m$, $ e_l\in \mathcal{E}_i^{out}$  and $\lambda_{i,0}\in \mathbf{R}^m$, and has a local variable $B_{i,0}=A_ix_{i,0}-b_{i,0}$ and a Poisson clock with rata $\varsigma_i$.

\noindent\rule{0.49\textwidth}{0.7mm}
{\bf Iteration at  $k$:}
Suppose player $i_k$'s clock ticks at time $k$, then player $i_k$ is active and updates its local
variables as follows:
\vskip 1mm
\quad {\bf Reading phase}: Get information from its neighbours' output buffer and its input buffer.
Duplicate its input buffer to its output buffer.
If $ j \in \mathbb{N}_{i_k}$, read ${x}^{(j)}_{k-\delta^k_j}$ and $\lambda_{j,k-\delta^k_j}$ from the output buffer of $j$.
If $e_l\in \mathcal{E}^{out}_{i_k}$ and $j \in \mathcal{N}_{l} \backslash i_k$,
then read $B_{j,k-\delta^k_{j}}$ and $z_{q,k-\pi^k_q}$, $e_q \in \mathcal{E}_{j}$ from the output buffer of $j$.

\vskip 1mm

\quad {\bf Computing phase}:

\qquad i): update $\lambda_{i_k}$,
\begin{eqnarray}
\tilde{\lambda}_{i_k,k} & = & \lambda_{i_k,k-\delta^k_{i_k}}+\sigma( A_{i_k}x_{i_k,k-\delta^k_{i_k}}-b_{i_k} \nonumber \\
\quad & \quad & + \sum_{e_l \in \mathcal{E}_{i_k}}V_{i_kl}z_{l,k-\pi^k_l})\label{equ_delayed_alg_pdi_1}\\
\lambda_{i_k,k+1}&=&\lambda_{i_k,k-\delta^k_{i_k}}+\eta(\tilde{\lambda}_{i_k,k}-\lambda_{i_k,k-\delta^k_i})\label{equ_delayed_alg_pdi_2}
\end{eqnarray}

\qquad ii): update $x_{i_k}$ and ${x}^{(i_k)}_{-i_k}$,
\begin{eqnarray}
\tilde{x}_{i_k,k}             &=&  P_{\Omega_{i_k}}\Big[ x_{i_k,k-\delta^k_{i_k}}
 -\tau\big( A_{i_k}^T(2\tilde{\lambda}_{i_k,k}-\lambda_{i_k,k-\delta^k_{i_k}}) \nonumber \\
\quad & \quad & +\nabla_{i_k}f_{i_k}(x_{i_k,k-\delta^k_{i_k}},{x}^{(i_k)}_{-i_k, {k-\delta^k_{i_k}}})  \nonumber \\
\quad & \quad & + \sum_{j\in \mathbb{N}_{i_k}}(x_{i_k,k-\delta^k_{i_k}}- {x}^{(j)}_{i, k-\delta^k_{j}} )  \big)  \Big]\label{equ_delayed_alg_pdi_x_1} \\
x_{i_k,k+1}                   &=&  x_{i_k,k-\delta^k_{i_k}}+\eta(\tilde{x}_{i_k,k}-x_{i_k,k-\delta^k_{i_k}})\label{equ_delayed_alg_pdi_x_2}\\
{x}^{(i_k)}_{-i_k,k+1} &=&  {x}^{(i_k)}_{-i_k,k-\delta^k_{i_k}} \nonumber\\
\quad & \quad & - \eta \tau \sum_{j\in \mathbb{N}_{i_k}}({x}^{(i_k)}_{-i_k,k-\delta^k_{i_k}}- {x}^{(j)}_{-i_k, k-\delta^k_{j}})\label{equ_delayed_alg_pdi_x_3}
\end{eqnarray}

\qquad iii): If $e_l\in \mathcal{E}^{out}_{i_k}$ and  $ j_k \in \mathbb{N}_{i_k}$, $j_k \in \mathcal{N}_{l} \backslash i_k$,
then update $z_{l,k}$:
\begin{eqnarray}
{z}_{l,k+1}        =   z_{l,k-\pi^k_{l}}- \eta\gamma(\lambda_{j_k,k-\delta^k_{j_k}}-\lambda_{i_k,k-\delta^k_{i_k}})\nonumber \\
-2\eta\sigma\gamma\big[(B_{j_k,k-\delta^j_{j_k}}-B_{i_k,k-\delta^k_{i_k}})+\sum_{q\in \mathcal{N}^e_{l}}L^e_{lq}z_{q,k-\pi^k_{q}}\big]\label{equ_delayed_alg_pdi_z_3}
\end{eqnarray}

\vskip 1mm
\quad {\bf Writing phase}:
 Write $x_{i_k,k+1}$, ${x}^{(i_k)}_{-i_k,k+1}$,  $\lambda_{i_k,k+1}$, $z_{l,k+1},e_l
 \in \mathcal{E}^{out}_{i_k}$  and $B_{i_k,k+1}=A_{i_k}x_{i_k,k+1}-b_{i_k}$ to its output buffer, and write $z_{l,k+1}$, $e_l\in \mathcal{E}^{out}_{i_k}$ to the  input buffer of player $j_k \in \mathcal{N}_{l} \backslash i_k$.

\vskip 3mm
 All other players keep their variables (output buffers) unchanged, and only increase $k$ to $k+1$.
\noindent\rule{0.49\textwidth}{0.7mm}
In \eqref{equ_delayed_alg_pdi_x_3}, ${x}^{(j)}_{-i_k, k-\delta^k_{j}}$ is player $j$'s estimation of all player's decision except $i_k$, including its own decision $x_j$.
All other notations are similar as the ones in  Algorithm \ref{alg_d}.
\end{alg}

\begin{rem}
The update for $\lambda_i$ and $z_l$ is the same for Algorithm \ref{alg_d} and Algorithm \ref{alg_d_partial information}.
So let us just state the differences between Algorithm \ref{alg_d} and Algorithm \ref{alg_d_partial information}.
\begin{itemize}
\item Algorithm \ref{alg_d} requires $i_k$ to read all players' decisions that  its  $f_{i_k}$ directly depends on, i.e.,  ${x}_{-i_k}$,
meanwhile Algorithm \ref{alg_d_partial information} only requires  player $i_k$ to inquire its neighbors' estimation ${x}^{(j)}, j\in \mathbb{N}_{i_k}$.
\item Let compare the update rule for local decision $x_{i_k}$, i.e., \eqref{equ_delayed_alg_2} and \eqref{equ_delayed_alg_pdi_x_1}-\eqref{equ_delayed_alg_pdi_x_2}.
\eqref{equ_delayed_alg_pdi_x_1} uses its local estimation ${x}^{(i_k)}_{-i_k}$ and local decision $x_{i_k}$ to calculate the gradient $\nabla_{i_k} f_{i_k}$
while \eqref{equ_delayed_alg_2} has to
use the true (delayed)  decisions of all other players ${x}_{-i_k}$. But \eqref{equ_delayed_alg_pdi_x_1} has an additional consensus term
$\sum_{j\in \mathbb{N}_{i_k}}(x_{i_k}- {x}^{(j)}_{i} )$ that will drive its local decision towards its neighbors' estimation.
\item Algorithm \ref{alg_d_partial information} has an additional equation \eqref{equ_delayed_alg_pdi_x_3} based on consensus dynamics to update $i_k$'s estimation of other players' decision ${x}^{(i_k)}_{-i_k}$.
\end{itemize}
Later on we will show that with the help of local consensus terms in   \eqref{equ_delayed_alg_pdi_x_1} and \eqref{equ_delayed_alg_pdi_x_3}, all players' estimation ${x}^{(i)}$  will reach consensus and converge to the true GNE ${x}^*$ of game \eqref{GM}.
Algorithm \ref{alg_d_partial information} is {\it fully distributed} since each player only needs local data and local communication to perform the updates.

\end{rem}

\subsection{Algorithm development}\label{sub_sec_pdi_develop}

Similar as the analysis for Algorithm \ref{alg_d} in subsection \ref{sub_sec_expand},
we can show that the asynchronous Algorithm \ref{alg_d_partial information} can be  treated as  a randomized block-coordinate fixed point iteration with delayed information, while the fixed point iteration is a preconditioned forward-backward method with properly chosen operators.
Let us make this clear step by step.

\vskip 5mm
Firstly, the synchronous version of Algorithm \ref{alg_d_partial information} without any delays for  agent $i$ at time $k$ is given as:
\begin{equation}
\begin{array}{l}\label{equ_syn_no_delay_pdi}
\tilde{\lambda}_{i,k}  =  \lambda_{i,k}+\sigma( A_{i}x_{i,k}  + \sum_{e_l \in \mathcal{E}_{i}}V_{il}z_{l,k}-b_{i}) \\
\lambda_{i,k+1}=\lambda_{i,k}+\eta(\tilde{\lambda}_{i,k}-\lambda_{i,k})\\
\tilde{x}_{i,k}       =  P_{\Omega_{i}}\Big[ x_{i,k}-\tau(\nabla_{i}f_{i}(x_{i,k},{x}^{(i)}_{-i, k})\\
+ A_i^T(2\tilde{\lambda}_{i,k}-\lambda_{i,k}) + \sum_{j\in \mathbb{N}_{i}} (x_{i,k}- {x}^{(j)}_{i,k} ) )  \Big]\\
x_{i,k+1}       =  x_{i,k}+\eta(\tilde{x}_{i,k}-x_{i,k})\\
{x}^{(i)}_{-i,k+1} =  {x}^{(i)}_{-i,k}- \eta\tau \sum_{j\in \mathbb{N}_{i}} ({x}^{(i)}_{-i,k}- {x}^{(j)}_{-i,k} ) )\\
{z}_{l,k+1}        =   z_{l,k}- \eta\gamma(\lambda_{j,k}-\lambda_{i,k})-2\eta\sigma\gamma(B_{j,k}-B_{i,k})\\
\quad -2 \eta\sigma\gamma\sum_{q\in \mathcal{N}^e_{l}}L^e_{lq}z_{q,k}
\end{array}
\end{equation}

{ To write all the players' update equations like \eqref{equ_syn_no_delay_pdi} in a compact form, we need some compact notations.
The compact notations  ${x}_{k}$,
 $\tilde{{x}}_k$,
 ${\Lambda}_k$,
 ${\tilde{\Lambda}}_k$,
 $\mathbf{b}$,
$\mathbf{z}_k$,
$\tilde{\mathbf{z}}_k$,
$\mathbf{A}$, $\mathbf{V}$ are defined the same as in subsection \ref{sub_sec_expand}.
Denote ${\mathbf{x}}_k= col({x}^{(1)}_k, {x}^{(2)}_k,\cdots, {x}^{(N)}_k) $ as the vector of all players' estimations on the overall decision profile,
and denote ${\mathbf{x}}_{-,k}=col( {x}^{(1)}_{-1,k}, {x}^{(2)}_{-2,k},\cdots, {x}^{(N)}_{-N,k}   )$.
Similar as $\tilde{\mathbf{z}}_k$,  we also  introduce intermediate variables
$ {\tilde{\mathbf{x}}}_k = col(\tilde{{x}}^{(1)}_k, \tilde{{x}}^{(2)}_k,\cdots, \tilde{{x}}^{(N)}_k) $ and
${\tilde{\mathbf{x}}}_{-,k}=col( \tilde{{x}}^{(1)}_{-1,k}, \tilde{{x}}^{(2)}_{-2,k},\cdots, \tilde{{x}}^{(N)}_{-N,k}   )$.
Denote $\mathbf{\mathbf{L}}= L\otimes I_n$.}

And then, we  define
\begin{equation} \label{equ_extend_pseudo_gradient}
\mathbf{F}({\mathbf{x}})=col(\nabla_{x_1} f_1(x_1,{x}^{(1)}_{-1}),\cdots, \nabla_{x_N} f_N(x_N,{x}^{(N)}_{-N})).
\end{equation}
Clearly,  $\mathbf{F}({\mathbf{x}})$ is different from  the pseudo-gradient $F({x})$ in \eqref{pseudogradient},
 since
$F({x})$ is evaluated at the true decision profile ${x}$ of all players while $\mathbf{F}({\mathbf{x}})$ is
evaluated at the estimations of the overall decision profile for each player.  Moreover, $\mathbf{F}({\mathbf{x}})$ is a mapping from $\mathbf{R}^{nN}$ to $\mathbf{R}^{n}$.
Therefore, we call $\mathbf{F}({\mathbf{x}})${\it extended pseudo-gradient }of game \eqref{GM}.

Motivated by \cite{pavel5}, we define two selection matrices $\mathcal{R}_i, \; \mathcal{S}_i$ as follows
\begin{eqnarray}
\mathcal{R}_i &=& [\begin{array}{ccc} \mathbf{0}_{n_i\times n_{<i}} & I_{n_i} & \mathbf{0}_{n_i\times n_{>i}}\end{array} ]\\
\mathcal{S}_i &=& \Big [ \begin{array}{ccc}
                  I_{n_{<i}} & \mathbf{0}_{n_{<i}\times n_i} & \mathbf{0}_{n_{<i}\times n_{>i}} \\
                  \mathbf{0}_{n_{>i}\times n_{<i}} & \mathbf{0}_{n_{>i}\times n_i} & I_{n_{>i}}
                \end{array} \Big ]
\end{eqnarray}
where $n_{<i}= \sum_{j<i,j\in \mathcal{N}}n_j$ and $n_{>i}=\sum_{j>i,j\in \mathcal{N}}n_j$.
Therefore, $\mathcal{R}_i {x}$ will  select $x_i$ from ${x}$ and $\mathcal{S}_i {x}$ will select ${x}_{-i}$ from ${x}$.
We stack $\mathcal{R}_i$ and $\mathcal{S}_i$ like $\mathbf{A}$ and denote $\mathcal{R}=\diag\{\mathcal{R}_1,\cdots,\mathcal{R}_N\}$, $S= \diag\{\mathcal{S}_1,\cdots,\mathcal{S}_N\}$.
Clearly,
 $$
{x}_k=\mathcal{R}{\mathbf{x}}_k,\;
{\mathbf{x}}_{-,k}=\mathcal{S}{\mathbf{x}}_k, \;
\tilde{{x}}_{k}=\mathcal{R} {\tilde{\mathbf{x}}}_k,\;
{\tilde{\mathbf{x}}}_{-,k}=\mathcal{S} {\tilde{\mathbf{x}}}_k.
$$

Then  with the above notations, all the players' update equations for the synchronous version of Algorithm \ref{alg_d_partial information}  can be written in a compact form as follows.
\begin{eqnarray}
{\tilde{\Lambda}}_{k}      &=& {\Lambda}_{k}+\sigma(\mathbf{A} \mathcal{R}{\mathbf{x}}_{k} +  \mathbf{V} \mathbf{z}_k-\mathbf{b}), \label{equ_syn_compact_pdi_1}\\
\tilde{\mathbf{z}}_{k}         &= &\mathbf{z}_{k}-\gamma \mathbf{V}^T {\Lambda}_{k} - 2\sigma\gamma\mathbf{V}^T (\mathbf{A}  \mathcal{R}{\mathbf{x}}_{k}-\mathbf{b})\nonumber\\
 \qquad &\qquad&-2 \sigma\gamma \mathbf{L}^e \mathbf{z}_k,\label{equ_syn_compact_pdi_2} \\
\tilde{{x}}_{k}         &=&P_{\Omega}\big[   \mathcal{R}{\mathbf{x}}_k-\tau(\mathbf{F}({\mathbf{x}}_k)\nonumber\\
 \qquad &\qquad& +\mathbf{A}^T(2{\tilde{\Lambda}}_{k}-{\Lambda}_{k}) + \mathcal{R}{\mathbf{L}}{\mathbf{x}}_k    )\big],\label{equ_syn_compact_pdi_3}\\
{\tilde{\mathbf{x}}}_{-,k} &=&  {\mathbf{x}}_{-,k}-\tau \mathcal{S}{\mathbf{L}}{\mathbf{x}}_k\label{equ_syn_compact_pdi_4}
\end{eqnarray}
and
\begin{equation}
\begin{array}{l}
{{\Lambda}}_{k+1}          =   {\Lambda}_k  +  \eta({\tilde{\Lambda}}_{k}-{\Lambda}_k ),  \quad
{\mathbf{z}}_{k+1}             =   \mathbf{z}_{k}   +  \eta (\tilde{\mathbf{z}}_{k}- \mathbf{z}_{k}),      \\
{{\mathbf{x}}}_{k+1}       =   {\mathbf{x}}_k     +  \eta ({\tilde{\mathbf{x}}}_{k}- {\mathbf{x}}_k  ). \label{equ_compact_relation_pdi}
\end{array}
\end{equation}

{ \begin{rem}
From the compact form \eqref{equ_syn_compact_pdi_1}-\eqref{equ_syn_compact_pdi_4}, we notice that Algorithm \ref{alg_d_partial information} utilizes three
graph-related matric:  incidence matrix $V$, Laplacian $L$ and edge Laplacian $L^e$. They play different roles in the algorithm.
The  Laplacian $L$ helps the players' local estimations to reach consensus on the GNE decision, while $V$ and $L^e$ help to ensure the coupling constraint and the consensus of local multipliers.
\end{rem}
}

\vskip 5mm
Then we can show that the compact from \eqref{equ_syn_compact_pdi_1}-\eqref{equ_syn_compact_pdi_4} and \eqref{equ_compact_relation_pdi}  can be derived as the fixed point iteration of
\begin{equation}\label{equ_fixed_pdi}
\mathbf{W}_{k+1} = \mathbf{W}_{k}+\eta (\mathcal{T} \mathbf{W}_k- \mathbf{W}_k ).
\end{equation}
where $\mathbf{W}$ denotes the stacked vector $\mathbf{W}=col( {\Lambda},\mathbf{z}, {\mathbf{x}}   )$.
Here,
\begin{equation}
{\bar{\mathfrak{A}}}: \mathbf{W} \mapsto
\left(
  \begin{array}{c}
    \mathbf{0}    \\
    \mathbf{0}  \\
    \mathcal{R}^TN_{\bar{\Omega}}(\mathcal{R}{\mathbf{x}}) \\
  \end{array}
\right)+
\left(
  \begin{array}{ccc}
    \mathbf{0} & -\mathbf{V} & -\mathbf{A} \mathcal{R} \\
    \mathbf{V}^T & \mathbf{0} & \mathbf{0} \\
    \mathcal{R}^T \mathbf{A}^T & \mathbf{0} & \mathbf{0}\\
  \end{array}
\right)\mathbf{W}
,\label{equ_operator_A_pdi}
\end{equation}

\begin{equation}
{\bar{\mathfrak{B}}}: \mathbf{W} \mapsto
\left(
  \begin{array}{c}
 \mathbf{b}    \\
   \mathbf{0 } \\
   \mathcal{R}^T \mathbf{F}({\mathbf{x}})+{\mathbf{L}}{\mathbf{x}} \\
  \end{array}
\right),\label{equ_operator_B_pdi}
\end{equation}

\begin{equation}
\bar{\Phi}=\left(
       \begin{array}{ccc}
         \sigma^{-1}I & \mathbf{V}       & \mathbf{A}\mathcal{R} \\
         \mathbf{V}^T    & \gamma^{-1}I  & \mathbf{0} \\
         \mathcal{R }^T\mathbf{A}^T & \mathbf{0}    & \tau^{-1}I \\
       \end{array}
     \right),\label{equ_operator_phi_pdi}
     \end{equation}
     \begin{equation}
\mathcal{T}= ({\rm Id}-\bar{\Phi}^{-1}\bar{\mathfrak{A}} )^{-1}({\rm Id}-\bar{\Phi}^{-1}\bar{\mathfrak{B}}).\label{equ_operator_T_pdi}
\end{equation}

Suppose $\bar{\Phi}$ is positive definite with properly chosen step-sizes.
Denote $\tilde{\mathbf{W}}=col({\tilde{\Lambda}}_k, \tilde{\mathbf{z}}_k,{\tilde{\mathbf{x}}}_k)$, then we only need to show that
$\tilde{\mathbf{W}}_k=\mathcal{T} \mathbf{W}_k=({\rm Id}-\bar{\Phi}^{-1}\bar{\mathfrak{A}} )^{-1}({\rm Id}-\bar{\Phi}^{-1}\bar{\mathfrak{B}})\mathbf{W}_k$ is equivalent with
 \eqref{equ_syn_compact_pdi_1}-\eqref{equ_syn_compact_pdi_4}. In other words, we need to show that
\begin{equation}\label{equ_forward_backward_splitting_pdi}
-\bar{\mathfrak{B}}\mathbf{W}_{k}\in \bar{\mathfrak{A}}\tilde{\mathbf{W}}_{k} +\bar{\Phi}(\tilde{\mathbf{W}}_{k}-\mathbf{W}_k).
\end{equation}
leads to \eqref{equ_syn_compact_pdi_1}-\eqref{equ_syn_compact_pdi_4}.

The derivation of  \eqref{equ_syn_compact_pdi_1} for  ${\Lambda}$ and \eqref{equ_syn_compact_pdi_2} for $\mathbf{z}$ is similar to subsection \ref{sub_sec_expand}.
The first component in  \eqref{equ_forward_backward_splitting_pdi} gives the update of ${\Lambda}$ as
$$
-\mathbf{b}= -\mathbf{A} \mathcal{R}{\tilde{\mathbf{x}}}_k-\mathbf{V}\tilde{\mathbf{z}}_k + \sigma^{-1}({\tilde{\Lambda}}_k-{\Lambda}_k)+\mathbf{V}(\tilde{\mathbf{z}}_k-\mathbf{z}_k)+\mathbf{A} \mathcal{R} ({\tilde{\mathbf{x}}}_k-{\mathbf{x}}_k), \nonumber
$$
which yields \eqref{equ_syn_compact_pdi_1}.
The $\mathbf{z}$ component in \eqref{equ_forward_backward_splitting_pdi} reads as
$$
\mathbf{0}= \mathbf{V}^T{\tilde{\Lambda}} + \mathbf{V}^T({\tilde{\Lambda}}_k-{\Lambda}_k )+ \gamma^{-1}({\tilde{\mathbf{z}}}_k-{\mathbf{z}}_k),  \nonumber
$$
hence,
$$
{\tilde{\mathbf{z}}}_k={\mathbf{z}}_k-\gamma\mathbf{V}^T(2{\tilde{\Lambda}}_k-{\Lambda}_k). \nonumber
$$
From  \eqref{equ_syn_compact_pdi_1} we have $2{\tilde{\Lambda}}_k-{\Lambda}_k= {\Lambda}_{k}+2\sigma(\mathbf{A} \mathcal{R }{\mathbf{x}}_{k} +  \mathbf{V} \mathbf{z}_k-\mathbf{b})$.
Then,
$$
{\tilde{\mathbf{z}}}_k={\mathbf{z}}_k-\gamma\mathbf{V}^T({\Lambda}_{k}+2\sigma(\mathbf{A} \mathcal{R}{\mathbf{x}}_{k} +  \mathbf{V} \mathbf{z}_k-\mathbf{b})), \nonumber
$$
which yields \eqref{equ_syn_compact_pdi_2} by  using $L^e=V^TV$.

The third component of \eqref{equ_forward_backward_splitting_pdi} reads as
\begin{equation}
\begin{array}{l} \label{equ_expand_x_pdi}
 -\mathcal{R}^T \mathbf{F}({\mathbf{x}}_k) -  {\mathbf{L}}{\mathbf{x}}_k
\in \mathcal{R}^TN_{{\Omega}}(\mathcal{R}{\tilde{\mathbf{x}}}) \\
\quad +\mathcal{R}^T\mathbf{A}^T {\tilde{\Lambda}}_k +  \mathcal{R}^T\mathbf{A}^T ({\tilde{\Lambda}}_k-{\Lambda}_k)+\tau^{-1}  ({\tilde{\mathbf{x}}}_k-{\mathbf{x}}_k)
\end{array}
\end{equation}
Hence, there exists $\mathbf{v}\in N_{{\Omega}}(\mathcal{R}{\tilde{\mathbf{x}}})=N_{{\Omega}}( \tilde{{x}}) $ such that the above equation becomes an equality.
It has been shown in \cite{pavel5} that both $\mathcal{R}$ and $\mathcal{S}$ are full row rank matrices, and
\begin{eqnarray}
\mathcal{R}^T\mathcal{R}+\mathcal{S}^T\mathcal{S}=I_{Nn} \\
\mathcal{R}\mathcal{S}^T=\mathbf{0}_{n\times n }, \quad \mathcal{S}\mathcal{R}^T=\mathbf{0}_{(N-1)n \times (N-1)n}\\
\mathcal{R}\mathcal{R}^T=I_n, \quad \mathcal{S}\mathcal{S}^T=I_{(N-1)n}
\end{eqnarray}

Multiplying \eqref{equ_expand_x_pdi} by $\mathcal{R}$ and by $\mathcal{R}\mathcal{R}^T=I_n$,  $\mathcal{R}{\tilde{\mathbf{x}}}_k=\tilde{{x}}$, $\mathcal{R}{\mathbf{x}}_k={x}_k$, we have
\begin{equation}
\begin{array}{l} \label{equ_expand_x_pdi_2}
 - \mathbf{F}({\mathbf{x}}_k) -  \mathcal{R} {\mathbf{L}}{\mathbf{x}}_k
= \mathbf{v}  +\mathbf{A}^T (2{\tilde{\Lambda}}_k-{\Lambda}_k)+\tau^{-1} (\tilde{\mathbf{x}}_k-\mathbf{x}_k)
\end{array}
\end{equation}
$\mathcal{R} {\mathbf{L}}{\mathbf{x}}_k$ is just the consensus terms for ${x}_k$.  Since $\mathbf{v}\in N_{{\Omega}}(\tilde{{x}}_k) $,
$\tau \mathbf{v} \in N_{{\Omega}}(\tilde{{x}}_k) $.
  Since $({\rm Id}+N_{{\Omega}})^{-1}({x})=P_{\Omega}({x})$, we have
$$\tilde{{x}}_{k}         =P_{\Omega}\big[   \mathcal{R}{\mathbf{x}}_k-\tau(\mathbf{F}({\mathbf{x}}_k)+\mathbf{A}^T(2{\tilde{\Lambda}}_{k}-{\Lambda}_{k}) + \mathcal{R}{\mathbf{L}}{\mathbf{x}}_k    )\big],
$$ which is \eqref{equ_syn_compact_pdi_3} for updating  the decision profile ${x}$.

Multiplying \eqref{equ_expand_x_pdi} by $\mathcal{S}$ and by $\mathcal{S}\mathcal{R}^T=\mathbf{0}$,
${\mathbf{x}}_{-,k}=\mathcal{S}{\mathbf{x}}_k $.
${\tilde{\mathbf{x}}}_{-,k}=\mathcal{S} {\tilde{\mathbf{x}}}_k$, we have
\begin{equation}
\begin{array}{l} \label{equ_expand_x_pdi_2}
-  \mathcal{S} {\mathbf{L}}{\mathbf{x}}_k=\tau^{-1} \mathcal{S}  ({\tilde{\mathbf{x}}}_k-{\mathbf{x}}_k),
\end{array}
\end{equation}
which is just \eqref{equ_syn_compact_pdi_4} for updating the estimations ${\mathbf{x}}_{-,k}$.

\vskip 3mm

Therefore, the synchronous version of Algorithm \ref{alg_d_partial information} can be regarded as the fixed point iteration \eqref{equ_fixed_pdi}, which finds fixed points of $\mathcal{T}$, or equivalently,
the zeros of $\bar{\mathfrak{A}}+\bar{\mathfrak{B}}$ when $\bar{\Phi}$ is positive definite.
Then we have the following  result characterizing  zeros of $\bar{\mathfrak{A}}+\bar{\mathfrak{B}}$.

\begin{thm}\label{thm_limiting_point_ine_1_pdi}
Suppose that Assumption \ref{assum1} holds for game \eqref{GM}, and assume the VI in \eqref{vi_original} has nonempty finite solution.
Then, any zero $\mathbf{W}^*= col({\Lambda}^*,\mathbf{z}^*,{\mathbf{x}}^*)$ of $\bar{\mathfrak{A}}+\bar{\mathfrak{B}}$  has the ${\mathbf{x}}^*=\mathbf{1}_N\otimes {x}^*$, and ${x}^*$ is  a variational GNE of \eqref{GM}, and ${\Lambda}^*=\mathbf{1}_N\otimes \lambda^*$, and ${x}^*$ and $\lambda^*$ satisfy KKT condition \eqref{kkt2} for VI \eqref{vi_original}.
\end{thm}

{\bf Proof:}
Suppose $\mathbf{W}^*= col({\Lambda}^*,\mathbf{z}^*,{\mathbf{x}}^*)$  is zero of  $\bar{\mathfrak{A}}+\bar{\mathfrak{B}}$,
then there exists $\mathbf{v}^*\in N_{{\Omega}}(\mathcal{R}{\mathbf{x}}^*)$ such that
\begin{eqnarray}
&&\mathbf{0} = -\mathbf{V}\mathbf{z}^*-\mathbf{A} \mathcal{R} {\mathbf{x}}^* + \mathbf{b}  \label{equ_zero_pdi_lambda}\\
&&\mathbf{0} = \mathbf{V}^T{\Lambda}^*  \label{equ_zero_pdi_z} \\
&&\mathbf{0} = \mathcal{R}^T [\mathbf{v}^* + \mathbf{A}^T {\Lambda}^* +\mathbf{F}({\mathbf{x}}^*)]+ {\mathbf{L}}{\mathbf{x}}^*  \label{equ_zero_pdi_x}
\end{eqnarray}
Since the communication graph $\mathcal{G}$ is  connected, $ \mathbf{V}^T{\Lambda}^* =\mathbf{0} $ of \eqref{equ_zero_pdi_z} implies that ${\Lambda}^*=\mathbf{1}_N\otimes \lambda^*$. We have $\mathbf{1}^T_N\otimes I_n \mathcal{R}^T= I_n$ and $\mathbf{1}^T_N\otimes I_n {\mathbf{L}}=\mathbf{0}$, hence, multiplying \eqref{equ_zero_pdi_x} with $\mathbf{1}^T_N\otimes I_n$ we have
\begin{equation}
\mathbf{v}^* + \mathbf{A}^T {\Lambda}^* +\mathbf{F}({\mathbf{x}}^*)=\mathbf{0}
\end{equation}
which in turn implies that ${\mathbf{L}}{\mathbf{x}}^*=\mathbf{0}$ with \eqref{equ_zero_pdi_x}.
Therefore, $ {\mathbf{x}}^*=\mathbf{1}_N\otimes {x}^*$ since $\mathcal{G}$ is connected.
Therefore, $\mathcal{R}{\mathbf{x}}^*= {x}^*$,  $\mathbf{v}^*\in N_{{\Omega}}({x}^*) $, and
$\mathbf{F}({\mathbf{x}}^*)=\mathbf{F}(\mathbf{1}_N\otimes {x}^*)=F({x}^*)$.
Then multiplying \eqref{equ_zero_pdi_lambda} with $\mathbf{1}_N\otimes I_m$, and multiplying \eqref{equ_zero_pdi_x}
with $\mathcal{R}$ we have
\begin{eqnarray}
&&\sum_{i=1}^N A_ix_i^*=\sum_{i=1}^N b_i\\
&&\mathbf{0} \in  N_{\Omega_i}(x_i^*) + A^T_i\lambda^*+ \nabla_{x_i} f_i(x_i^*, {x}^*_{-i}), i\in \mathcal{N}.
\end{eqnarray}
which is exactly the KKT condition \eqref{kkt2} for the variational GEN of game \eqref{GM}.
\hfill $\Box$

\subsection{Algorithm convergence analysis}\label{sub_sec_pdi_convergence}

With the algorithm development analysis in subsection\ref{sub_sec_pdi_develop},
we can define another random variable $\xi_k$ similar as subsection \ref{sec_sub_har} such that
Algorithm \ref{alg_d_partial information} can be written as
\begin{equation}\label{equ_fixed_pdi_asy}
\mathbf{W}_{k+1} = \mathbf{W}_{k}+\eta \xi_k \odot (\mathcal{T} \hat{\mathbf{W}}_k- \mathbf{W}_k ).
\end{equation}
where $\hat{\mathbf{W}}_k$ is a vector with delayed information.
Therefore,  to show the convergence of Algorithm \ref{alg_d_partial information} we only need to show that the randomized block-coordinate fixed iteration with delayed information \eqref{equ_fixed_pdi_asy} is convergent.
To give a sufficient convergence condition,
we impose the following assumption on the extended pseudo-gradient $ \mathbf{F}({\mathbf{x}})$ of game \eqref{GM}.

\begin{assum}\label{assum4}
The extended pseudo-gradient $ \mathbf{F}({\mathbf{x}})$ in \eqref{equ_extend_pseudo_gradient} has $\mathcal{R}^T\mathbf{F}({\mathbf{x}})$ to be $\bar{\chi}-$ cocoercive, i.e.,$\forall {\mathbf{x}}, {\mathbf{x}}^{'}  \in \mathbf{R}^{nN}$
\begin{equation}
\langle {\mathbf{x}} - {\mathbf{x}}^{'}, \mathcal{R}^T\mathbf{F}({\mathbf{x}})-\mathcal{R}^T\mathbf{F}({\mathbf{x}}^{'})\rangle \geq \bar{\chi} || \mathbf{F}({\mathbf{x}})-\mathbf{F}({\mathbf{x}}^{'}) ||_2^2,
\end{equation}
\end{assum}

\begin{rem}
In fact, Assumption \ref{assum4} is not equivalent with Assumption \ref{assum2}. So the class of problems that Algorithm \ref{alg_d} works for might not be exactly the same
as the class of problems that Algorithm \ref{alg_d_partial information} works for.
 Assumption \ref{assum4} has been adopted in \cite{pavel5}, \cite{pavel6} and  \cite{tatarenko} for distributed NE/GNE computation, and it can be checked with the help of Jacobian condition.
\end{rem}

Then similar to Lemma \ref{lem_operator_A_B},  we have

\begin{lem}\label{lem_operator_A_B_pdi}
Suppose Assumption \ref{assum1} and \ref{assum4} hold. Denote $d^*=\max \{|L|_{11},\cdots,|L|_{NN}\}$ as the maximal degree of communication graph $\mathcal{G}$.
Operator $\bar{\mathfrak{A}}$ in \eqref{equ_operator_A_pdi} is maximally monotone. Operator $\bar{\mathfrak{B}}$ in \eqref{equ_operator_B_pdi} is $\beta_{E}$-cocoercive for
$\beta_{E}\in (0, \frac{1}{2} \min\{ \bar{\chi}, \frac{1}{2d^*}\})$.

\end{lem}

{\bf Proof:}
In fact, ${\bar{\mathfrak{A}}}$ is a sum of $\bar{\mathfrak{A}}_1: \mathbf{W} \mapsto  \mathbf{0} \times \mathbf{0}\times \mathcal{R}^T N_{{\Omega}}(\mathcal{R}{\mathbf{x}}) $ and a matrix operate  $\bar{\mathfrak{A}}_2$.  $\bar{\mathfrak{A}}_2$ is maximally monotone since the matric is skew-symmetric.
Since $ N_{{\Omega}}$ is maximally monotone, and matrix $\mathcal{R} $ is full row rank, the composition $\mathcal{R}^T \circ N_{{\Omega}} \circ \mathcal{R}$
is also maximally monotone.
In fact, we take $\mathbf{v} = col(v_1,\cdots,v_n)$ with $v_i \in N_{\Omega_i}(x_i)$, and $\mathbf{v}^{'}=col(v^{'}_1,\cdots,v^{'}_n)$ with $v^{'}_i\in N_{\Omega_i}(x^{'}_i)$
then $ \langle  {\mathbf{x}}-{\mathbf{x}}^{'}, \mathcal{R}^T (\mathbf{v}-\mathbf{v}^{'})\rangle= \sum_{i=1}^N \langle x_i-x^{'}_i, v_i-v^{'}_i  \rangle \geq 0$.
Hence, $\mathcal{R}^T \circ N_{{\Omega}} \circ \mathcal{R}$ is monotone, and
next we show that it  is maximally monotone. Suppose $({\mathbf{x}},\mathcal{R}^T\mathbf{u}) \in gra  \mathcal{R}^T \circ N_{{\Omega}} \circ \mathcal{R} $.
Then $\mathbf{u}=col(u_1,\cdots,u_n)$ with $u_i \in N_{\Omega_i}(x_i)$. Suppose we have another vector $({\mathbf{x}}^{'}, \mathcal{R}^T\mathbf{u}^{'})$ such that
$ \langle {\mathbf{x}}-{\mathbf{x}}^{'}, \mathcal{R}^T(\mathbf{u}- \mathbf{u}^{'}) \rangle \geq 0$, then we need to show that
$({\mathbf{x}}^{'}, \mathcal{R}^T\mathbf{u}^{'})\in gra \mathcal{R}^T \circ N_{{\Omega}} \circ \mathcal{R}$.
With $ \langle {\mathbf{x}}-{\mathbf{x}}^{'}, \mathcal{R}^T(\mathbf{u}- \mathbf{u}^{'}) \rangle \geq 0$ we know that
$\sum_{i=1}^N \langle x_i-x^{'}_i, u_i-u^{'}_i  \rangle \geq 0$. Since $u_i\in N_{\Omega_i}(x_i)$ and $N_{\Omega_i}$ is maximally monotone, then $u_i^{'}\in N_{\Omega_i}(x_i^{'})$. Since $\mathcal{R}$ is full row rank, $\mathbf{u}^{'}\in N_{\Omega}(\mathcal{R}{\mathbf{x}}^{'})$.
Therefore, ${\bar{\mathfrak{A}}}$ is maximally monotone.

For operator  ${\bar{\mathfrak{B}}}$, it can be written as ${\bar{\mathfrak{B}}}={\bar{\mathfrak{B}}}_1+{\bar{\mathfrak{B}}}_2$ with
${\bar{\mathfrak{B}}}_1(\mathbf{W}) = \mathbf{b} \times \mathbf{0 }\times    \mathcal{R}^T \mathbf{F}({\mathbf{x}})$ and
$ {\bar{\mathfrak{B}}}_2(\mathbf{W}) = \mathbf{0}\times \mathbf{0 }\times{\mathbf{L}}{\mathbf{x}} $.
With Assumption \ref{assum4}, we have
$
\langle \mathbf{W}-\mathbf{W}^{'}, {\bar{\mathfrak{B}}}_1(\mathbf{W})-{\bar{\mathfrak{B}}}_1(\mathbf{W}^{'}) \rangle
=
\langle  \bar{\mathbf{x}}-\bar{\mathbf{x}}^{'}, \mathcal{R}^T\mathbf{F}({\mathbf{x}})-\mathcal{R}^T\mathbf{F}({\mathbf{x}}^{'})\rangle
\geq
 \bar{\chi} || \mathbf{F}({\mathbf{x}})-\mathbf{F}({\mathbf{x}}^{'}) ||_2^2
= \bar{\chi} || \mathcal{R}^T\mathbf{F}({\mathbf{x}})-\mathcal{R}^T\mathbf{F}({\mathbf{x}}^{'}) ||_2^2
=
\bar{\chi} || {\bar{\mathfrak{B}}}_1(\mathbf{W})-{\bar{\mathfrak{B}}}_1(\mathbf{W}^{'}) ||_2^2
$ where we use the fact that $\mathcal{R}\mathcal{R}^T=I_n$.

Since ${\mathbf{L}}$ is semi-positive definite, ${\mathbf{L}}{\mathbf{x}}$ can be regarded as the gradient of $\frac{1}{2}{\mathbf{x}}^T{\mathbf{L}}{\mathbf{x}}$.
Therefore, we have
$  \langle  {\mathbf{x}}-{\mathbf{x}}^{'}, {\mathbf{L}}{\mathbf{x}}-{\mathbf{L}}{\mathbf{x}}^{'}  \rangle \geq
\frac{1}{|| {\mathbf{L}}||_2} ||  {\mathbf{L}}{\mathbf{x}}-{\mathbf{L}}{\mathbf{x}}^{'}||_2^2$ with Bailon-Haddad theorem.
Therefore,
$
\langle \mathbf{W}-\mathbf{W}^{'}, {\bar{\mathfrak{B}}}_2(\mathbf{W})-{\bar{\mathfrak{B}}}_2(\mathbf{W}^{'}) \rangle
 =
 \langle{\mathbf{x}}-{\mathbf{x}}^{'}, {\mathbf{L}}{\mathbf{x}}-{\mathbf{L}}{\mathbf{x}}^{'}  \rangle
 \geq
\frac{1}{|| {\mathbf{L}}||_2} || {\mathbf{L}}{\mathbf{x}}-{\mathbf{L}}{\mathbf{x}}^{'}||_2^2
\geq
\frac{1}{2d^*} || {\bar{\mathfrak{B}}}_2(\mathbf{W})-{\bar{\mathfrak{B}}}_2(\mathbf{W}^{'}) ||_2^2
$ where we use the fact that $||{\mathbf{L}} ||_2 \leq 2d^* $.

Therefore, with $|| a||_2^2+ ||b ||_2^2 \geq \frac{1}{2} ||a+b||_2^2 $,
we have
 $
\langle \mathbf{W}-\mathbf{W}^{'}, {\bar{\mathfrak{B}}}(\mathbf{W})-{\bar{\mathfrak{B}}}(\mathbf{W}^{'}) \rangle
\geq
 \bar{\chi} || {\bar{\mathfrak{B}}}_1(\mathbf{W})-{\bar{\mathfrak{B}}}_1(\mathbf{W}^{'}) ||_2^2
 +
\frac{1}{2d^*} || {\bar{\mathfrak{B}}}_2(\mathbf{W})-{\bar{\mathfrak{B}}}_2(\mathbf{W}^{'}) ||_2^2
\geq
\frac{1}{2} \min\{ \bar{\chi}, \frac{1}{2d^*}\} || {\bar{\mathfrak{B}}}(\mathbf{W})-{\bar{\mathfrak{B}}}(\mathbf{W}^{'}) ||_2^2
$. Therefore, $\bar{\mathfrak{B}}$ in \eqref{equ_operator_B_pdi} is $\beta_{E}$-cocoercive for
$\beta_{E}\in (0, \frac{1}{2} \min\{ \bar{\chi}, \frac{1}{2d^*}\})$.
\hfill $\Box$

Then similar to Lemma \ref{lem_preconditioned_operators}, we have,
\begin{lem}\label{lem_preconditioned_operators_pdi}
Suppose Assumption\ref{assum1}, \ref{assum4} hold.
Take $ \delta > \frac{1}{2\beta_{E}}$ where ${\beta_{E}}\in  (0, \frac{1}{2} \min\{ \bar{\chi}, \frac{1}{2d^*}\}) $, and take the step-sizes $\tau, \gamma, \sigma$ in Algorithm \ref{alg_d_partial information}  such that $\bar{\Phi}-\delta I$ is positive semi-definite for $\bar{\Phi}$ in \eqref{equ_operator_phi_pdi}.
Then under the $\bar{\Phi}-$induced norm $||\cdot||_{\bar{\Phi}}$, we have,
(i): $\bar{\Phi}^{-1}\bar{\mathfrak{A}}$ is  maximally monotone, and $\mathcal{T}_1=({\rm Id}+\bar{\Phi}^{-1}\bar{\mathfrak{A}})^{-1}\in \mathcal{A}(\frac{1}{2})$.
(ii): $\bar{\Phi}^{-1}\bar{\mathfrak{B}}$ is $ \delta\beta_{E}-$cocoercive, and  $\mathcal{T}_2={\rm Id}-\bar{\Phi}^{-1}\bar{\mathfrak{B}}\in \mathcal{A}(\frac{1}{2\delta\beta_{E}})$.
(iii): $\mathcal{T}$ in \eqref{equ_operator_T_pdi} is  an averaged operator, and $\mathcal{T} \in \mathcal{A}(\frac{2\delta\beta_{E}}{4\delta\beta_{E}-1})  $.
\end{lem}
{\bf Proof:} The proof of (i) and (ii) follow from Lemma 5.6 of \cite{yipeng}, and  (iii) follows from Proposition 2.4 of \cite{combettes3}.
\hfill $\Box$

\vskip 3mm
Then combined with Assumption \ref{assum3},  we obtain a sufficient condition for the convergence of the asynchronous GNE computation algorithm with partial decision information, i.e. Algorithm \ref{alg_d_partial information}.
\begin{cor}\label{cor_1}
Suppose Assumption \ref{assum1}, \ref{assum2}, \ref{assum3} and \ref{assum4} hold.
Take $ \delta > \frac{1}{2\beta_{E}}$ where ${\beta_{E}}\in  (0, \frac{1}{2} \min\{ \bar{\chi}, \frac{1}{2d^*}\})$, and take the step-sizes $\tau, \gamma, \sigma$ in Algorithm \ref{alg_d_partial information}  such that $\bar{\Phi}-\delta I$ is positive semi-definite for $\bar{\Phi}$ in \eqref{equ_operator_phi_pdi}.
Denote $p_{min}=\min_{i\in \mathcal{N}}\{\frac{\varsigma_i}{\sum_{j=1}^N \varsigma_j}\}$, and choose $0<\eta \leq  \frac{c N p_{\min}}{2\Psi\sqrt{p_{\min}}+1} \frac{4\delta\beta_{E}-1}{2\delta\beta_{E}}$ for any $c\in (0,1)$.
Then with Algorithm \ref{alg_d_partial information}, with probability $1$ the players' local estimations  $x^{(i)}_{k},i\in \mathcal{N}$ asymptotically reach consensus, i.e.,
$\lim_{k\rightarrow \infty} {x}^{(i)}_k={x}^*, \forall i\in \mathcal{N}$,
and the players' local multipliers ${\lambda}_{i,k},i\in \mathcal{N}$ reach consensus, i.e.,
$\lim_{k\rightarrow \infty} \lambda_{i,k}=\mathbf{\lambda}^*, \forall i\in \mathcal{N}$.
 Moreover, the limiting point ${x}^*$ is a variational GNE  of game in \eqref{GM}, which
  together with  $\lambda^*$ satisfy KKT \ref{kkt2}.
\end{cor}

{\bf Proof:}
The result can be regarded as a  corollary of  Theorem \ref{thm_convergence_asy}. With the algorithm development analysis in subsection \ref{sub_sec_pdi_develop}, the distributed asynchronous algorithm with partial decision information, i.e. Algorithm \ref{alg_d_partial information}, can be written in a compact form as \eqref{equ_fixed_pdi_asy}
with a delayed vector $\hat{W}_k$. With Assumption \ref{assum3}, the delays is bounded for all the time steps.
With Lemma \ref{lem_preconditioned_operators_pdi}, the operator $\mathcal{T}$ is  an averaged operator, and $\mathcal{T} \in \mathcal{A}(\frac{2\delta\beta_{E}}{4\delta\beta_{E}-1})  $.
Hence,  the almost surely convergence of \eqref{equ_fixed_pdi_asy} to a fixed point of $\mathcal{T}$
is guaranteed by Lemma 13 and Theorem  14 of \cite{yinwotao_arock} if the step-size $\eta$ is chosen accordingly.
Since the fixed point of $\mathcal{T}$ is the same as the zero of $\bar{\mathfrak{A}}+\bar{\mathfrak{B}}$, then the limiting points of Algorithm \ref{alg_d_partial information} enjoy the desirable properties with Theorem \ref{thm_limiting_point_ine_1_pdi}.
\hfill $\Box$

\section{Numerical studies}\label{sec_numerical_studies}

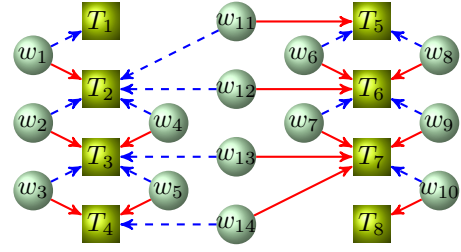
\begin{figure}
\begin{center}
\begin{tikzpicture}[->,>=stealth',shorten >=0.3pt,auto,node distance=2cm,thick,
  rect node/.style={rectangle, ball color={rgb:red,0;green,0.2;yellow,1},font=\sffamily,inner sep=1pt,outer sep=0pt,minimum size=14pt},
  wave/.style={decorate,decoration={snake,post length=0.1mm,amplitude=0.5mm,segment length=3mm},thick},
  main node/.style={shape=circle, ball color=green!20,text=black,inner sep=1pt,outer sep=0pt,minimum size=15pt},scale=0.9]


  \foreach \place/\i in {{
  (-3,0.5)/1},
  {(-3,-0.5)/2},
  {(-3,-1.5)/3},
  {(-1,-0.5)/4},
  {(-1,-1.5)/5},
  {(1,0.5)/6},
  {(1,-0.5)/7},
  {(3,0.5)/8},
  {(3,-0.5)/9},
  {(3,-1.5)/10},
  {(0,1)/11},
  {(0,0)/12},
  {(0,-1)/13},
  {(0,-2)/14}}
    \node[main node] (a\i) at \place {};

      \node at (-3,0.5) {\rm \color{black}{$w_{1}$}};
      \node at (-3,-0.5){\rm \color{black}{$w_2$}};
      \node at (-3,-1.5){\rm \color{black}{$w_3$}};
      \node at (-1,-0.5){\rm \color{black}{$w_{4}$}};
      \node at (-1,-1.5){\rm \color{black}{$w_{5}$}};
      \node at (1,0.5)  {\rm \color{black}{$w_{6}$}};
      \node at (1,-0.5) {\rm \color{black}{$w_7$}};
      \node at (3,0.5)  {\rm \color{black}{$w_8$}};
      \node at (3,-0.5) {\rm \color{black}{$w_9$}};
      \node at (3,-1.5) {\rm \color{black}{$w_{10}$}};
      \node at(0,1)     {\rm \color{black}{$w_{11}$}};
      \node at (0, 0)   {\rm \color{black}{$w_{12}$}};
      \node at (0,-1)   {\rm \color{black}{$w_{13}$}};
      \node at (0,-2)   {\rm \color{black}{$w_{14}$}};

  \foreach \place/\x in {{(-2,1)/1},{(-2,0)/2},{(-2,-1)/3},
    {(-2,-2)/4}, {(2,1)/5}, {(2,0)/6}, {(2,-1)/7}, {(2,-2)/8}}
  \node[rect node] (b\x) at \place {};

      \node at (-2,1) {\rm \color{black}{$T_1$}};
      \node at (-2,0){\rm \color{black}{$T_2$}};
      \node at (-2,-1)  {\rm \color{black}{$T_3$}};
      \node at (-2,-2) {\rm \color{black}{$T_4$}};
      \node at (2,1){\rm \color{black}{$T_5$}};
      \node at (2,0) {\rm \color{black}{$T_6$}};
      \node at (2,-1)  {\rm \color{black}{$T_7$}};
      \node at (2,-2) {\rm \color{black}{$T_8$}};

  \path[dashed,->,blue,thick]               (a1) edge (b1);
  \path[->,red,thick]                       (a1) edge (b2);

  \path[dashed,->,blue,thick]               (a2) edge (b2);
  \path[->,red,thick]                       (a2) edge (b3);

  \path[dashed,->,blue,thick]               (a3) edge (b3);
  \path[->,red,thick]                       (a3) edge (b4);

  \path[dashed,->,blue,thick]               (a4) edge (b2);
  \path[->,red,thick]                       (a4) edge (b3);

  \path[dashed,->,blue,thick]               (a5) edge (b3);
  \path[->,red,thick]                       (a5) edge (b4);

  \path[dashed,->,blue,thick]               (a6) edge (b5);
  \path[->,red,thick]                       (a6) edge (b6);

  \path[dashed,->,blue,thick]               (a7) edge (b6);
  \path[->,red,thick]                       (a7) edge (b7);

  \path[dashed,->,blue,thick]               (a8) edge (b5);
  \path[->,red,thick]                       (a8) edge (b6);

  \path[dashed,->,blue,thick]               (a9) edge (b6);
  \path[->,red,thick]                       (a9) edge (b7);

  \path[dashed,->,blue,thick]               (a10) edge (b7);
  \path[->,red,thick]                       (a10) edge (b8);

  \path[dashed,->,blue,thick]               (a11) edge (b2);
  \path[->,red,thick]                       (a11) edge (b5);

  \path[dashed,->,blue,thick]               (a12) edge (b2);
  \path[->,red,thick]                       (a12) edge (b6);

  \path[dashed,->,blue,thick]               (a13) edge (b3);
  \path[->,red,thick]                       (a13) edge (b7);

  \path[dashed,->,blue,thick]               (a14) edge (b4);
  \path[->,red,thick]                       (a14) edge (b7);
\end{tikzpicture}
\end{center}
\caption{Task allocation game: An edge from $w_i$ to $T_j$ on this graph implies that a part of  worker $w_i$'s output is  allocated to task $T_j$.}\label{fig_task_allocation_game}
\end{figure}

In this part, we consider a task allocation game with $8$ tasks $\{T_1,\cdots,T_8\}$ and $14$ processors (workers) $\{w_1,\cdots,w_{14}\}$ as an illustrative example
of the proposed asynchronous GNE seeking algorithms.
Each task $T_j$ is quantified as a load of $C_j>0$ that should be met by the workers.
Each worker $w_i$ decides its working output $x_i=col(x^1_{i},x^2_{i},x^3_{i}, x^4_{i})\in \mathbf{R}^4$ within its capacity $\mathbf{0} \leq x_i\leq B_i, B_i \in \mathbf{R}_{+}^4$.
If worker $w_i$ allocates a part of its output to task $T_j$, there is an arrow $w_i\rightarrow T_j$ on Fig. \ref{fig_task_allocation_game}.
Specifically, if $w_i$ allocates $x^1_{i},x^2_{i}$ to $T_j$, there is a dashed blue arrow on Fig. \ref{fig_task_allocation_game}, and
if $w_i$ allocates $x^3_{i},x^4_{i}$ to $T_j$, there is a solid red arrow on Fig. \ref{fig_task_allocation_game}.
Define a matrix $A=[A_1, \cdots, A_{15}]\in \mathbf{R}^{8\times 56}$
with $A_i=[a^1_i,a_i^2,a_i^3,a_i^4] \in \mathbf{R}^{8\times 4}$ quantifying how the output of worker $w_i$  is allocated to each task. Each column $a_i^k$ has only one element being nonzero, and the $j$th element of $a_i^1$ or $a_i^2$ is nonzero if there is a  dashed blue arrow $w_i\rightarrow T_j$ on Fig. \ref{fig_task_allocation_game},
and the $j$th element of $a_i^3$ or $a_1^4$ is nonzero if there is a  red solid arrow $w_i\rightarrow T_j$ on Fig. \ref{fig_task_allocation_game}.
The nonzero elements in $A_i$ are randomly chosen from $[0.5,1]$, which could be regarded as delivery loss factors.
It is required that
the tasks should be met by the working output of the players.
Denote $C=col(C_1,\cdots,C_8)$, then the workers have an equality coupling constraint:
$A{x}=C$ with $x=col(x_1,\cdots,x_{14})$.
The objective function of player (worker) $w_i$ is
\begin{equation}\label{equ_objective_example}
f_i(x_i,\mathbf{x}_{-i}) = c_i(x_i) - R^T(\mathbf{x}) A_i x_i.
\end{equation}
Here,  $c_i(x_i)$ is a cost function, which is  taken as
$c_i(x_i)=  \sum_{k=1}^4 q_i^k(x_i^k+1)\log(x_i^k+1)+(p_i^Tx_i-d_i)^2+x_i^T S_ix_i$.
$R({x})=col(R_1({x}),\cdots,R_8({x}))$
is a vector function that maps ${x}$ to the award price of each task, and
$R_j({x})=\kappa_j - \chi_j [A{x}]_j$.

The parameters of the game are generated as follows.
The following parameters are randomly drawn as
$C_j\in [1,2]$, $\chi_j \in [0.5,1]$, $\kappa_j \in [4,9]$, $q^k_i\in [0.5,1.5]$,  $d_i\in [1,2]$.
$p_i \in \mathbf{R}^{4}$ is a randomly generated stochastic vector,  $S_i\in \mathbf{R}^{4\times 4}$ is a randomly generated positive definite matrix,  and each element of $B_i$ is  drawn from $[1,2]$.
The parameters are numerically checked to ensure Assumption \ref{assum1} and \ref{assum2}, and to ensure additional Assumption \ref{assum4} for partial decision information case with the Jacobian condition in \cite{pang,scutari,combettes1}.

Then we give some common algorithm settings.
The players communicate over a ring graph Figure \ref{fig_multiplier_graph}, and edges are arbitrarily ordered to define the incidence matrix $V$ and edge Laplacian $L^e$.
We assume that $\varsigma_1=,\cdots,=\varsigma_N$, hence $\mathbb{P}(i_k=i)=\frac{1}{N}$ in both algorithms. The maximal delay $\Psi$ in Assumption \ref{assum3} is taken to to be a finite number while the delayed vector  $\hat{\varpi}_k$ in \eqref{equ_randomized_delayed} or $\hat{W}_k$ in \eqref{equ_fixed_pdi_asy} is constructed with the coordinates randomly chosen from the past $\Psi$ steps.
Each player has a local $C_i=\frac{1}{15}C$, and  the step-sizes $\sigma,\gamma,\tau$ are randomly chosen from $[0.3,0.5]$ such that $\Phi$ in \eqref{precon_phi} or $\bar{\Phi}$ in \eqref{equ_operator_phi_pdi} is positive definite,  and
$\eta=0.4$. The initial ${x}_{i,0}$ or the initial ${x}^{(i)}_{i,0}$ is randomly chosen within $[\mathbf{0} ,  B_i]$, and initial $\lambda_i$, $z_l$ are chosen to be zero. For the partial decision information case, the initial ${x}^{(i)}_{-i,0}$ is chosen to be $\mathbf{0}$.

\subsection{Simulation of the algorithm with full decision information}

\begin{figure}
\begin{center}
\begin{tikzpicture}[->,>=stealth',shorten >=0.3pt,auto,node distance=2.1cm,thick,
  rect node/.style={rectangle,ball color=blue!10,font=\sffamily,inner sep=1pt,outer sep=0pt,minimum size=12pt},
  wave/.style={decorate,decoration={snake,post length=0.1mm,amplitude=0.5mm,segment length=3mm},thick},
  main node/.style={shape=circle, ball color=green!20,,text=black,inner sep=1pt,outer sep=0pt,minimum size=12pt},scale=0.8]


  \foreach \place/\i in {
  {(-2,2)/1},
  {(-1,2)/2},
  {(0,2)/3},
  {(1,2)/4},
  {(2,2)/5},
  {(2,1)/6},
  {(2,0)/7},
  {(2,-1)/8},
  {(1,-1)/9},
  {(0,-1)/10},
  {(-1,-1)/11},
  {(-2,-1)/12},
  {(-2,0)/13},
  {(-2,1)/14}}
    \node[main node] (a\i) at \place {};

      \node at (-2,2){\rm \color{black}{$w_1$}};
      \node at (-1,2){\rm \color{black}{$w_2$}};
      \node at (0,2){\rm \color{black}{$w_3$}};
      \node at (1,2){\rm \color{black}{$w_4$}};
      \node at (2,2){\rm \color{black}{$w_5$}};
      \node at (2,1){\rm \color{black}{$w_6$}};
      \node at (2,0){\rm \color{black}{$w_7$}};
      \node at (2,-1){\rm \color{black}{$w_8$}};
      \node at (1,-1){\rm \color{black}{$w_9$}};
      \node at (0,-1){\rm \color{black}{$w_{10}$}};
      \node at (-1,-1){\rm \color{black}{$w_{11}$}};
      \node at (-2,-1){\rm \color{black}{$w_{12}$}};
      \node at (-2,0){\rm \color{black}{$w_{13}$}};
      \node at (-2,1) {\rm \color{black}{$w_{14}$}};

        \path[<->,blue,thick]               (a1) edge (a2);
         \path[<->,blue,thick]               (a2) edge (a3);
         \path[<->,blue,thick]               (a3) edge (a4);
         \path[<->,blue,thick]               (a4) edge (a5);

         \path[<->,blue,thick]               (a5) edge (a6);
         \path[<->,blue,thick]               (a6) edge (a7);
         \path[<->,blue,thick]               (a7) edge (a8);
         \path[<->,blue,thick]               (a8) edge (a9);
           \path[<->,blue,thick]               (a9) edge (a10);
         \path[<->,blue,thick]               (a10) edge (a11);
             \path[<->,blue,thick]             (a11) edge (a12);
            \path[<->,blue,thick]               (a12) edge (a13);
            \path[<->,blue,thick]               (a13) edge (a14);
          \path[<->,blue,thick]               (a1) edge (a14);

\end{tikzpicture}
\end{center}
\caption{Communication graph $\mathcal{G}$: Player $w_i$ and $w_j$ are able to exchange their local information if there exists an edge between them on this graph. }\label{fig_multiplier_graph}
\end{figure}
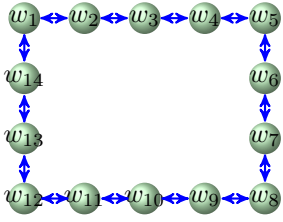

In this part, we examine the asynchronous algorithm with full decision information, i.e., Algorithm \ref{alg_d}.
The players exchange their local $\lambda_i$ and inquire  $z_l$ through the communication graph of Figure \ref{fig_multiplier_graph}.
The maximal delay  $\Psi=20$.
Moreover,   Algorithm \ref{alg_d} requires each player is able to know the decisions of the players who also contribute to the same task as itself.
For example, player  $w_2$ contributes to both task $T_2$ and task $T_3$ in Figure \ref{fig_task_allocation_game}. Meanwhile, $\{w_1,w_2,w_4, w_{11},w_{12} \}$ all contribute to Task $T_1$, and
$\{ w_2,w_3,w_4,w_{5},w_{13}\}$ all contribute to Task $T_2$. With the specific structure of the objective function \eqref{equ_objective_example},
player $w_2$ is required to be able to know the (delayed) decision of $\{w_1,w_3,w_4,w_{5}, w_{11},w_{12},w_{13}\}$.

\begin{figure}
  \centering
  \includegraphics[width=3.5in]{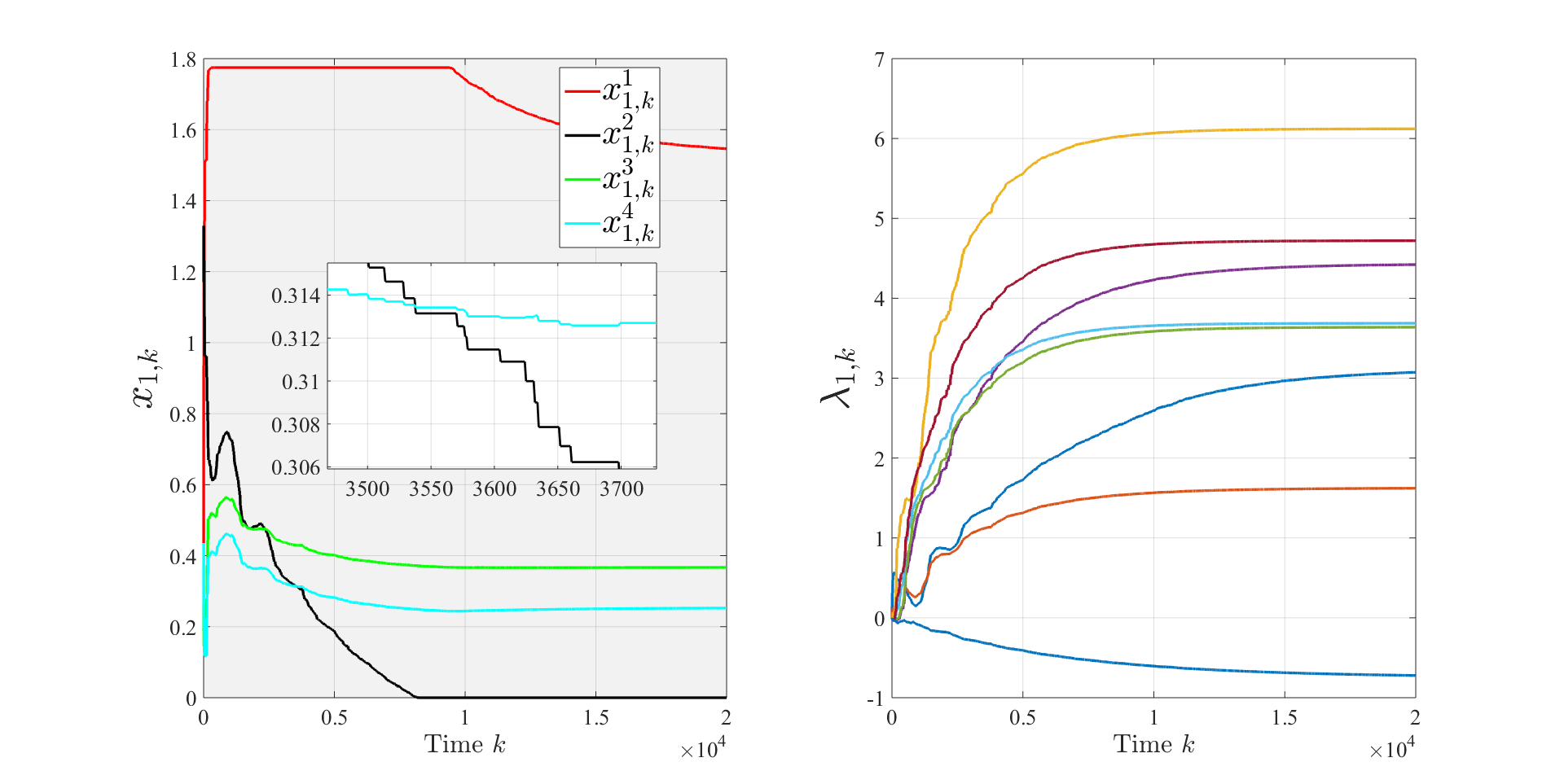}\\
  \caption{The trajectories of player $w_1$'s  $x_{1,k}$ and $\lambda_{1,k}$ generated by ADAGNES Algorithm \ref{alg_d}}\label{fig1}
\end{figure}

\begin{figure}
  \centering
  \includegraphics[width=3.5in]{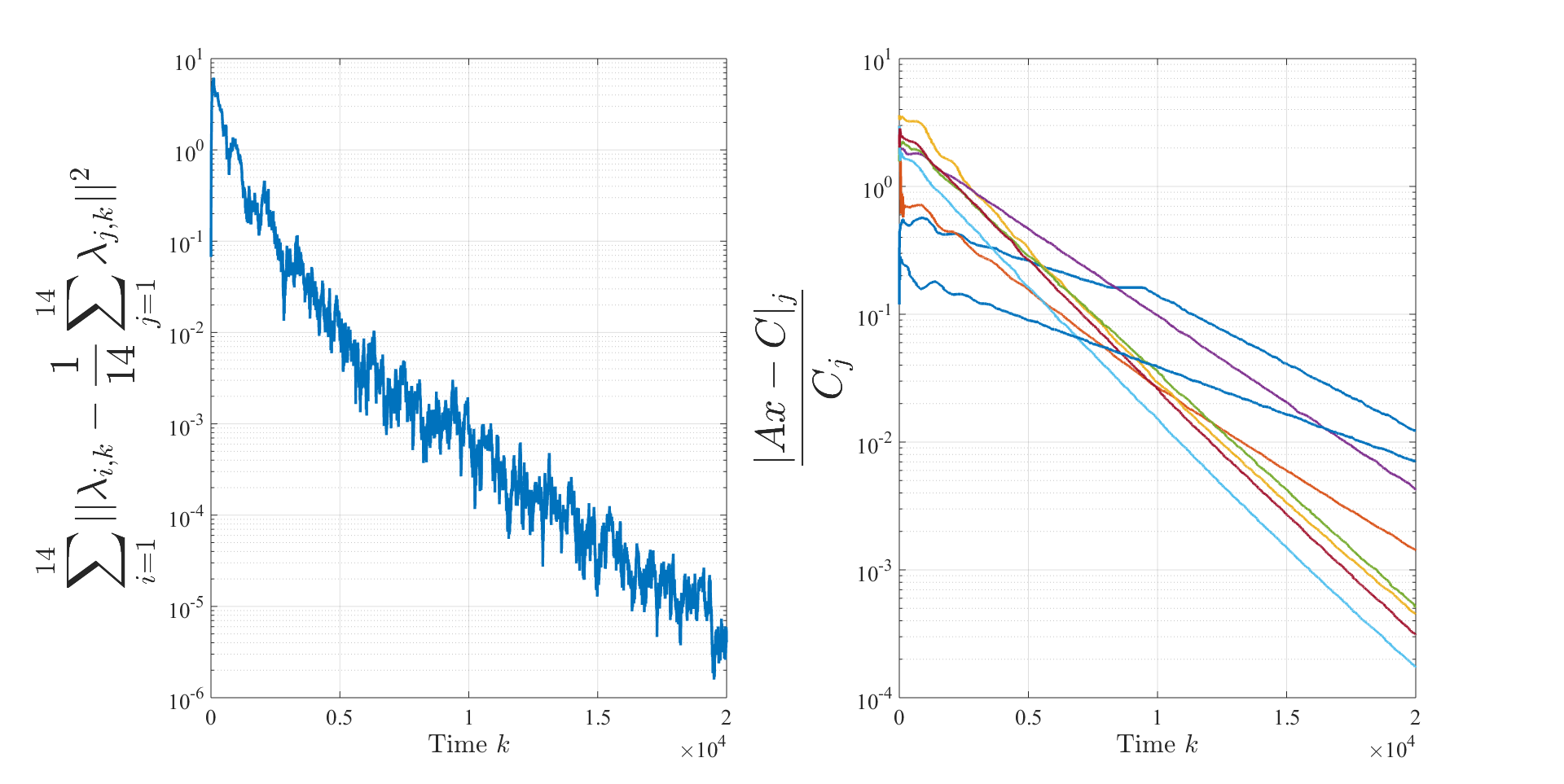}\\
 \caption{The trajectories of  the consensual error of all $\lambda_i$ and the violation of equality constraint generated by ADAGNES Algorithm \ref{alg_d}}\label{fig2}
\end{figure}

\begin{figure}
  \centering
  \includegraphics[width=3.5in]{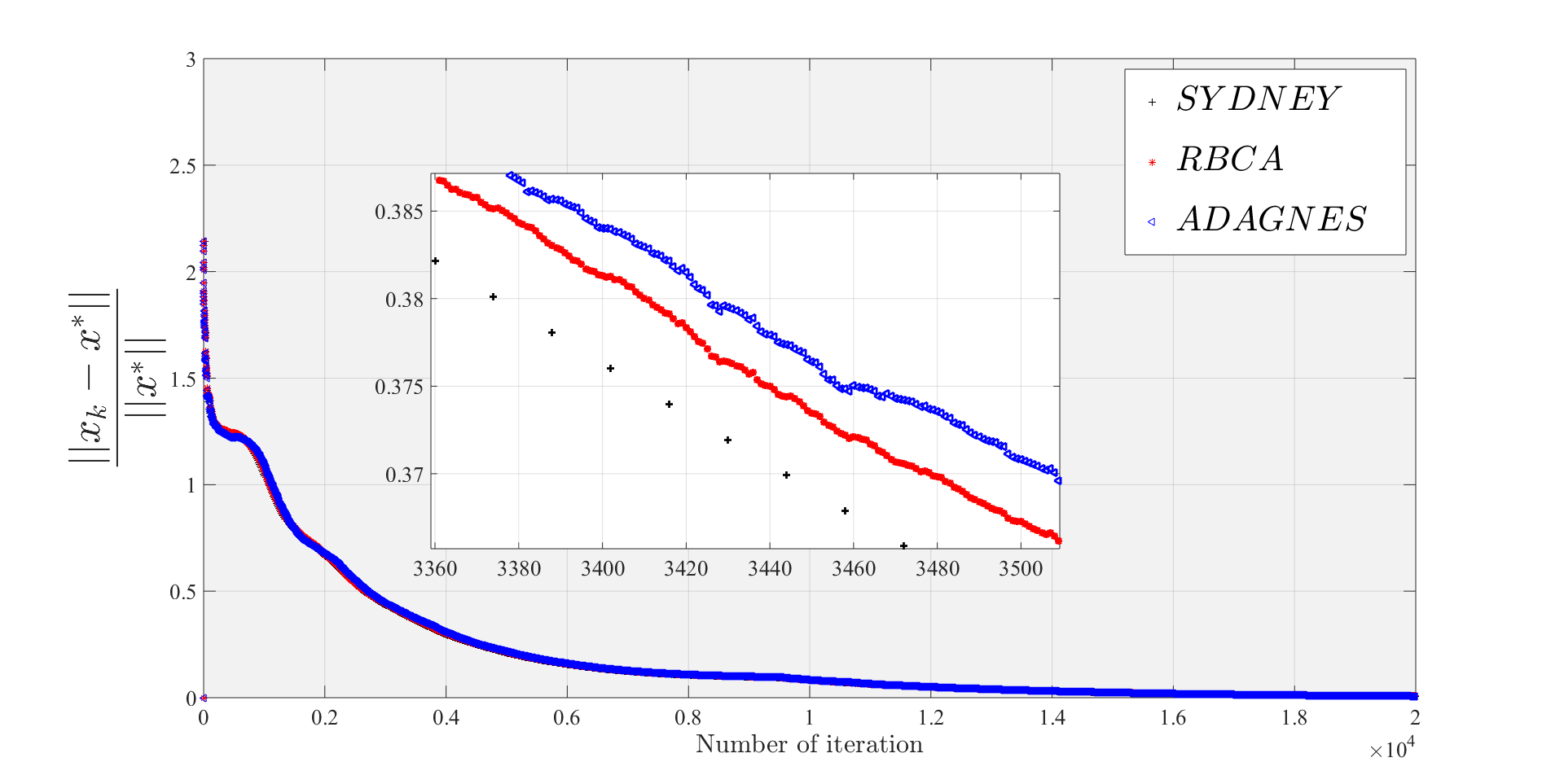}\\
   \caption{Comparison of algorithm efficiency in term of number of total iterations for SYDNEY, RBCA, and ADAGNES Algorithm \ref{alg_d}}\label{fig3}
\end{figure}

 We compare the performance of  the synchronous algorithm, i.e., SYDNEY, the randomized block-coordinated algorithm (RBCA), i.e.,  \eqref{equ_randomized}, and ADAGNES Algorithm \ref{alg_d},  by running them with the same parameters.
Fig \ref{fig1} and \ref{fig2} show trajectories generated by ADAGNES, demonstrating  its convergence.
However, since it is run on a single laptop, it takes a long iteration step for ADAGNES to  converge.
We can compare the algorithms'   performance by showing how $\frac{||{x}_k-{x}^* ||}{|| {x}^*||}$ is evolving with the number of all agents' total iterations. Note that in SYDNEY, all agents perform
iterations at each time step, hence at each time step, SYDNEY has N times, i.e, 14 times more computations and communications compared to RBCA and ADAGNES.  Figure \ref{fig3}shows $\frac{||{x}_k-{x}^* ||}{|| {x}^*||}$ vs the  number of all agents' total iterations for the three algorithms. It shows that  ADAGNES is only a bit of inferior to SYDNEY, and competitive with RBCA in terms of the number of agents' total iterations.
However, since ADAGNES removes  the players' idle time  and does not need any global coordination, it is expected that ADAGNES would be superior to SYDNEY and RBCA when implemented in distributed multi-agent systems.

\subsection{Simulation of the algorithm  with partial decision information}

In this part, we validate the asynchronous algorithm with partial decision information, i.e., Algorithm  \ref{alg_d_partial information}.
Each player can only communicate with its neighbors through the graph in Figure \ref{fig_multiplier_graph}.
For example, different from the previous subsection,  player $w_2$ cannot know the (delayed) decision of $\{w_1,w_3,w_4,w_{5}, w_{11},w_{12},w_{13}\}$, but can only
{ use its local estimation of all  players' decision profile $x^{(2)}$ and its neighbor $w_1,w_3$'s estimations $x^{(1)}$, $x^{(3)}$,  obtained by
communication over the graph in Figure \ref{fig_multiplier_graph}.}
The maximal delay is  $\Psi=30$. The other parameters are set the same as before.

\begin{figure}
  \centering
  \includegraphics[width=3.5in]{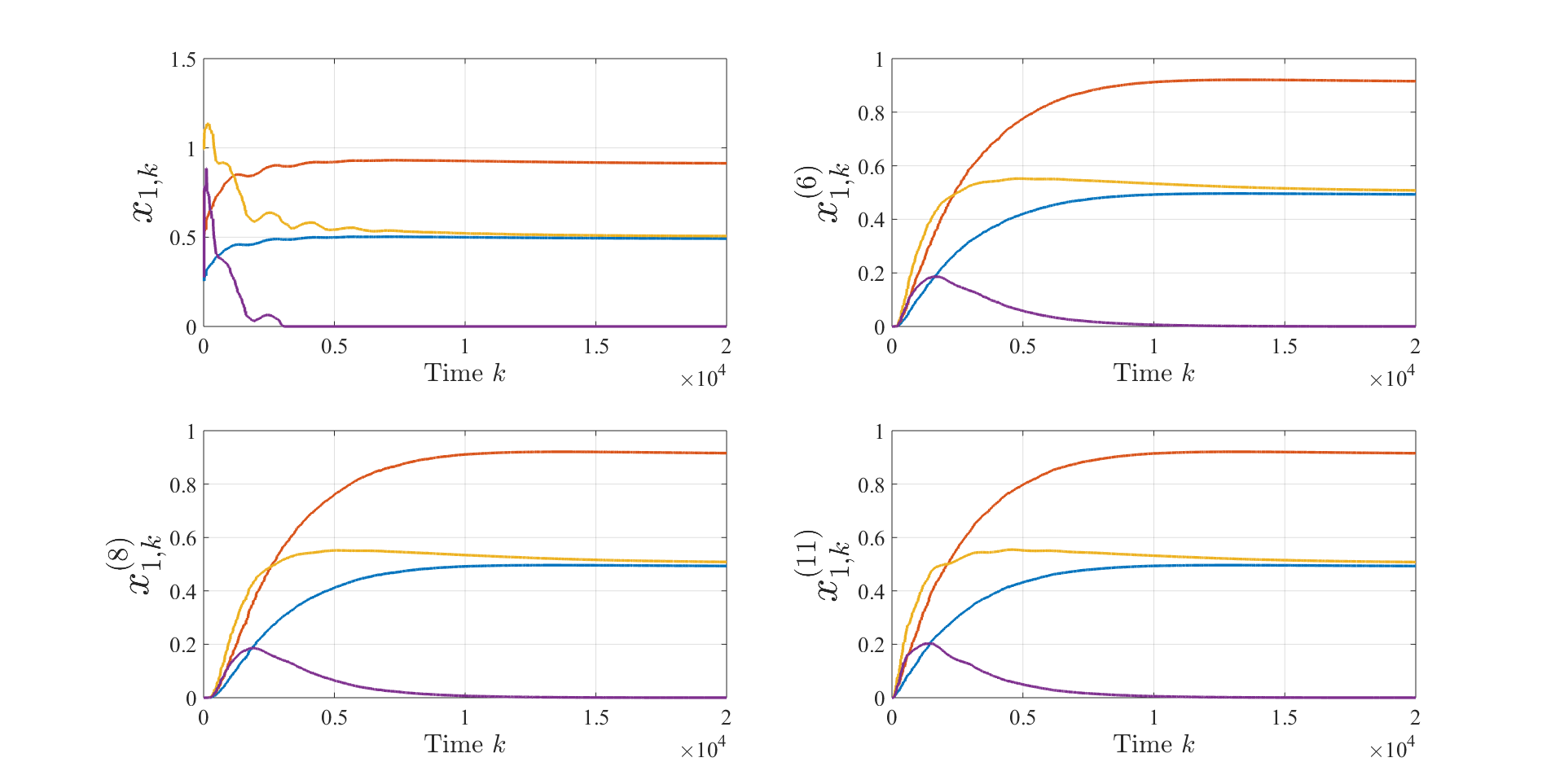}\\
  \caption{The trajectories of player $w_1$'s  decision $x_{1,k}$ and the trajectories of player $6$, $8$, $11$' estimation of player $1$'s decision, i.e.,  $x^{(6)}_{1,k}$, $x^{(8)}_{1,k}$ and $x^{(11)}_{1,k}$,  generated by ADAGNES-PDI Algorithm \ref{alg_d_partial information}}\label{fig_pdi_1}
\end{figure}

\begin{figure}
  \centering
  \includegraphics[width=3.5in]{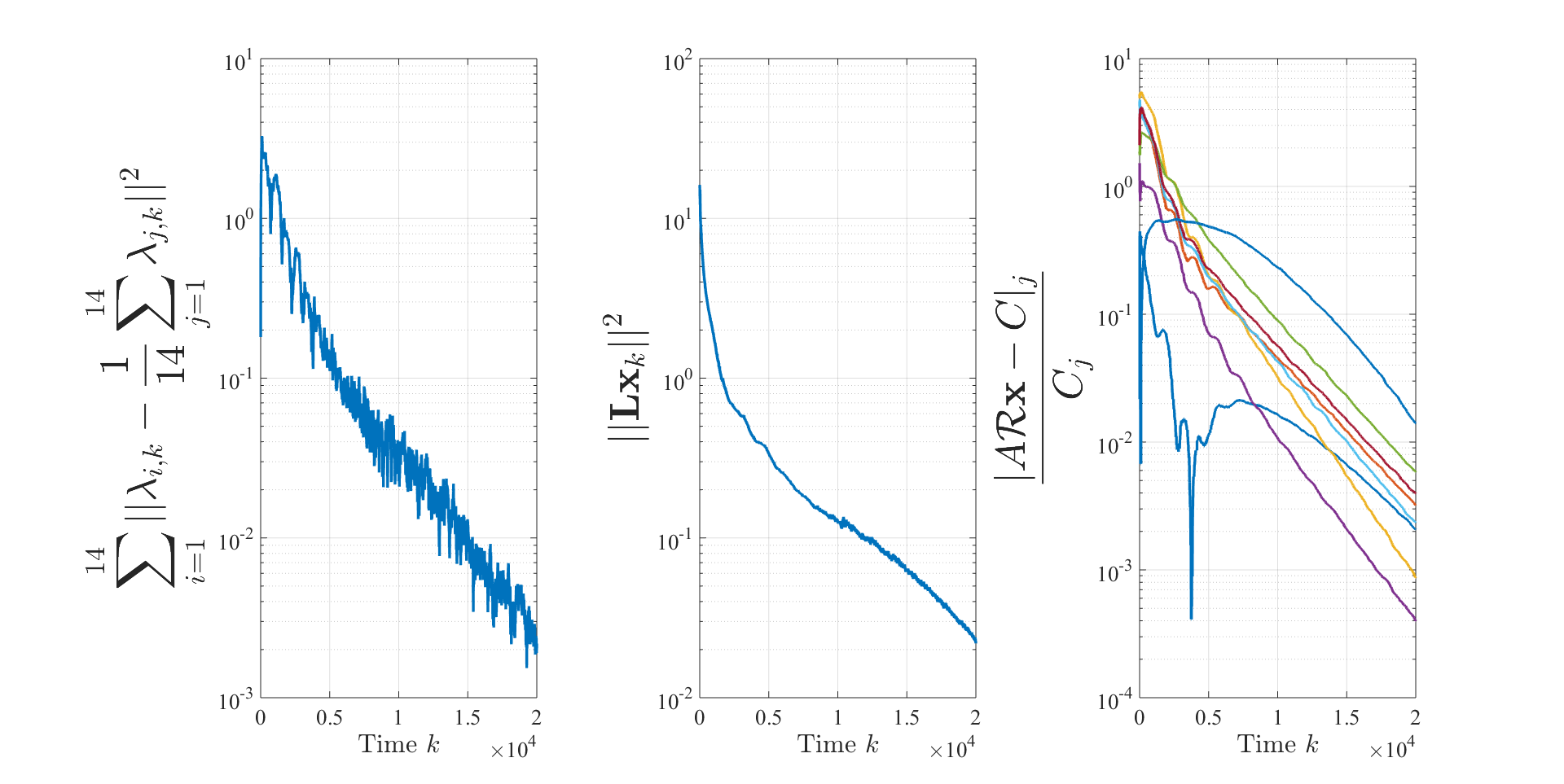}\\
  \caption{The trajectories of the  performance index  for consensus error of local multipliers, the  performance index  for consensus error of local decision estimations,
  and the violation of the affine constraints, all   generated by ADAGNES-PDI Algorithm \ref{alg_d_partial information}}\label{fig_pdi_2}
\end{figure}

Figure \ref{fig_pdi_1} shows the trajectories of player $w_1$' decision and the estimations of player $w_1$'s decision by player $w_6$, $w_8$ and $w_{11}$.
It demonstrates that for each player's decision, the other players' local estimation of it will asymptotically track the true decision of that player.
Figure \ref{fig_pdi_2} shows the trajectories of the consensus error of local multipliers, and the consensus error of the local decision estimations.
They verify that all the local multipliers  will reach consensus and all the local estimations for the overall decision profile will also reach consensus.
Figure \ref{fig_pdi_2} also shows the trajectories of the violation of the affine constraint by the true decision profile. It verifies that
the affine coupling constraint will be satisfied asymptotically.


\begin{figure}
  \centering
  \includegraphics[width=3.5in]{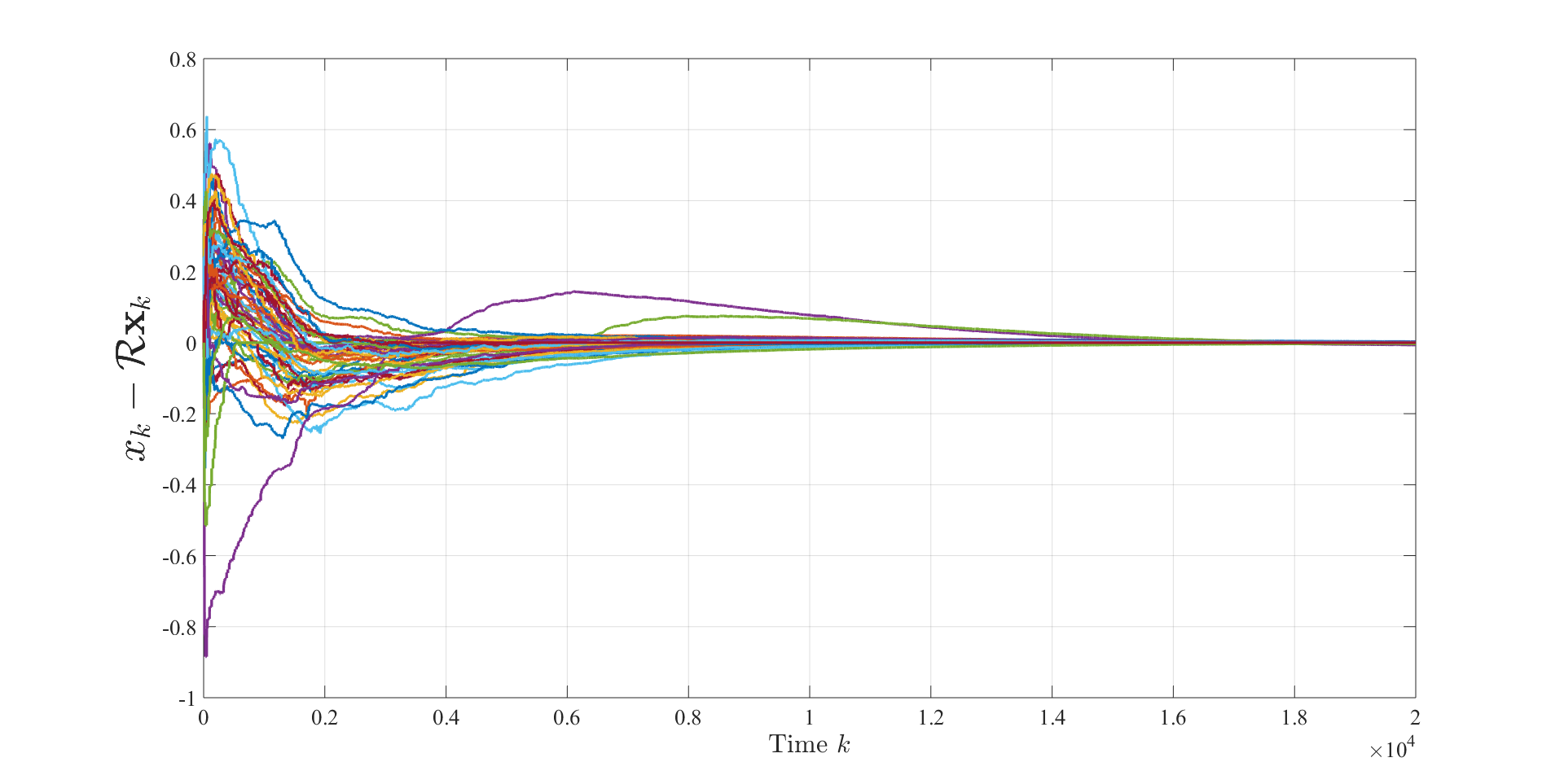}\\
  \caption{Comparison of the  Algorithm \ref{alg_d} and \ref{alg_d_partial information}:  ${x}_k$ is generated from Algorithm \ref{alg_d} and ${\mathbf{x}}_k$
  is generated from Algorithm \ref{alg_d_partial information}, with the same initial conditions and algorithm parameters. }\label{fig_pdi_3}
\end{figure}

To compare the asynchronous algorithms with full decision information and partial decision information, i.e., Algorithm \ref{alg_d} and Algorithm \ref{alg_d_partial information}, we run them with the same initial conditions and algorithm parameters.
Figure \ref{fig_pdi_3} shows the trajectories of the difference between the decision profile ${x}_k$ generated from Algorithm \ref{alg_d} and the decision profile $\mathcal{R}{\mathbf{x}}_k$  generated from Algorithm \ref{alg_d_partial information}.
It verifies that both algorithms have the overall decision profile  converge to same (unique) variational GNE.

\section{Conclusions}\label{sec_concluding}

In this paper, we proposed two asynchronous distributed  GNE seeking algorithms for noncooperative games with affine coupling constraints under full decision information and partial decision information, respectively.
In the full decision information case, with the help of edge variables and edge Laplacian matrix, the agent can perform local iteration with local data and delayed neighbour information without any synchronization or global coordination.
In the partial decision information case, the agent can use its local estimation to perform local gradient evaluation, while an consensus dynamics with Laplacian matrix ensures that the local estimations will reach the same GNE.
 Both algorithms are shown to be unified under the umbrella of preconditioned forward-backward methods for finding zeros of monotone operators.
 Then, sufficient fixed step-size choices are given for algorithm convergence by resorting to
  the theory of randomized fixed point iteration of non-expansive operators with delayed information. Our simulations verify only the convergence of the proposed algorithms. It would be interesting to perform the experiments on  actual distributed systems and study the efficiency of the asynchronous algorithm compared with synchronous ones. We also notice that the assumptions for the convergence analysis of the algorithm with {\it full decision information} and the algorithm with {\it partial decision information} are different. Hence, how to eliminate the gap should be a future research direction. Motivated by \cite{hu2}, it is also promising to investigate the asynchronous GNE computation with switching or directed communication graphs.

\end{document}